\RequirePackage{ifpdf}
\ifpdf 
\documentclass[pdftex]{sigma}
\else
\documentclass{sigma}
\fi

\usepackage{mathrsfs}
\usepackage{bbm}
\usepackage{stmaryrd}
\usepackage[all]{xy}

\newcommand{\im}{{\rm im}}
\DeclareMathOperator*{\cycl}{\textstyle\sum\hspace{-1.3em}\circlearrowleft}

\newcommand{\z}{\mathbb Z}
\renewcommand{\r}{\mathbb R}

\renewcommand{\c}{\mathbb C}

\newcommand{\X}{\mathfrak{X}}

\numberwithin{equation}{section}
\begin{document}
\allowdisplaybreaks

\renewcommand{\PaperNumber}{077}

\FirstPageHeading

\ShortArticleName{The Torsion of Spinor Connections and Related Structures} 

\ArticleName{The Torsion of Spinor Connections\\ and Related Structures}

\Author{Frank KLINKER} 
\AuthorNameForHeading{F. Klinker}
\Address{University of Dortmund, 44221 Dortmund, Germany}
\Email{\href{mailto:frank.klinker@math.uni-dortmund.de}{frank.klinker@math.uni-dortmund.de}}

\ArticleDates{Received August 25, 2006, in f\/inal form November
03, 2006; Published online November 09, 2006}

\Abstract{In this text we introduce the torsion of  spinor
connections. In terms of the tor\-sion we give conditions on a
spinor connection to produce Killing vector f\/ields. We relate
the Bianchi type identities for the torsion of spinor connections
with Jacobi identities for vector f\/ields on supermanifolds.
Furthermore, we discuss applications of this notion of torsion. }

\Classification{17B66; 53C27; 53B20}

\section{Introduction}
In this article we introduce the torsion of arbitrary spinor
connections. Although  the  construction depends on additional
data on the spinor bundle, namely a choice of charge conjugation,
the notion of torsion of a spinor connection is a natural
extension of what is usually known as the torsion of a connection
on a manifold. In Section~\ref{sectionTorsion} we give the
relevant def\/initions and discuss certain properties. In
particular, in Proposition~\ref{Bianchi} we list Bianchi-type
identities which connect the torsion and the curvature of the
given spinor connection.

\looseness=1 The spinor connections for which parallel spinors
leads to inf\/initesimal transformations of the underlying
manifold are discussed in Section~\ref{admissibleconnection}. This
turns out to be a symmetry condition on the torsion and lead to
the def\/inition of admissibility. In the case of metric
connections, admissibility recovers the connections with totally
skew symmetric torsion. The latter have been discussed in detail
during the last years, e.g.~\cite{IvanovIvanov} and references
therein. Beside these metric connections there are a lot of
examples coming from supergravity models and we  emphasize on
them, e.g.~\cite{CremmJulScher} for the basic one. In view of the
Fierz relation we formulate the admissibility condition in terms
of forms. In Theorem~\ref{admissibleForm} and its
extension~\ref{admissibleForm2}
 we give a list of all admissible connections,
 i.e.\ connections such that the supersymmetry bracket of parallel spinor f\/ields --
 when identif\/ied with the projection from the endomorphisms of the spinor
 bundle to the one-forms -- closes into the space of Killing vector f\/ields without
 further assumptions. Such connections are always used when we consider
 supergravity theories and exa\-mine  the variations of the odd f\/ields.
 Moreover admissible pairs are one of the basic objects in our current
 work  on natural  realizations of  supersymmetry  on non-f\/lat manifolds.
In  Section~\ref{JacobiBianchi} we draw a connection to the
geometry of a special class of supermanifold. We show that
torsion enters  naturally  into higher order commutators of
canonically def\/ined super vector f\/ields. This yields a
connection between the graded Jacobi identity on the superalgebra
of vector f\/ields and the Bianchi identities derived in
Section~\ref{sectionTorsion}. The motivation for the introduction
of supermanifolds and the consideration of the canonical vector
f\/ield is taken from the constructions in \cite{Klinker2} and
\cite{Klinker4}. The canonical vector f\/ield we consider has also
been discussed in \cite{Papado2} from another point  of view: One
of the vector f\/ields is considered as  f\/irst order
ope\-ra\-tor on the bundle of exterior powers of the spin bundle
and it is asked when this operator is a~dif\/ferential.

Section \ref{exmps} is devoted to examples and applications. We
introduce three  notions  of torsion freeness which are motivated
by the discussion so far, and we shortly discuss torsion freeness
in the case of f\/lat space. Some properties of spinor connections
on f\/lat space have recently been discussed in \cite{BernNagy}.
In Section \ref{examplebranes} we discuss  brane metrics admitting
torsion free admissible subsets.

\section{Preliminaries}\label{inclusions}

We consider the graded manifold $\hat M=(M,\Gamma\Lambda S)$.
where $M$ denotes a (pseudo) Riemannian spin manifold and $\Lambda
S$ the  exterior bundle of the spinor bundle $S$. The splitting
$\Lambda S=\Lambda_0S\oplus\Lambda_1 S$ into even and odd forms
def\/ine the even and odd functions on $\hat M$. An  inclusion of
vector f\/ields on the base manifold $M$ and sections of the
spinor bundle $S$ into the vector f\/ields on $\hat M$ via
$\jmath:\X(M)\oplus \Gamma S\hookrightarrow \X(\hat M)$ yields a
splitting\footnote{We often  use the identif\/ications
$\Gamma(E\oplus F) =\gamma E\oplus\Gamma F$, $\Gamma(E\otimes
F)=\Gamma E\otimes_{C^\infty(M)} \Gamma F$, $\Gamma({\rm
Hom}(E,F))={\rm Hom}_{C^\infty(M)}(\Gamma E,\Gamma F)$ etc, for
sections of vector bundles over the manifold $M$.}
\begin{equation*}
\X(\hat M) = \Gamma\Lambda  S\otimes \X(M)\oplus \Gamma\Lambda
S\otimes \Gamma S,
\end{equation*}
compare \cite{Kostant} or \cite{ManSardan}. The even and odd parts
of the vector f\/ields are given by
\[
\X(\hat M)_\alpha = \Gamma\Lambda_{\alpha\,{\text{mod}}2} S\otimes
\X(M)\oplus \Gamma\Lambda_{\alpha+1\,{\text{mod}}2}S\otimes \Gamma
S, \qquad\alpha=0,1.
\]
The v-  and  s-like f\/ields are def\/ined by
\[
\X_{\text{v}}(\hat M):=\Gamma\Lambda
S\otimes\X(M),\qquad\X_{\text{s}}(\hat M):=\Gamma\Lambda  S\otimes
\Gamma S.
\]
We call a vector f\/ield $X$ of order $(k,1)$ or  $(k,0)$ if
$X\in\Gamma\Lambda^k S\otimes \Gamma S$ or $ X\in\Gamma\Lambda^k
S\otimes \X(M)$, respectively.

The graded manifold $\hat M=(M,\Gamma\Lambda  S)$ is equipped with
a bilinear form $g+C$ where $g$ is the metric on $M$ and $C$ a
charge conjugation on $S$. The latter is a spin-invariant bilinear
form on~$S$. Another important map is the Clif\/ford
multiplication
\[
\gamma: \X(M)\otimes \Gamma S\to \Gamma S, \qquad
\gamma(X\otimes\eta)=\gamma(X)\eta = X\eta
\]
with
\begin{equation*}
XY+YX=-2g(X,Y).
\end{equation*}
As the notation indicates, we often consider the induced map
$\gamma:\X(M)\to \Gamma\, {\rm End}( S)$. We call the images of a
local frame $\{e_\mu\}$ on $M$  $\gamma$-matrices and write
$\gamma(e_\mu)=\gamma_\mu$. We always use the abbreviation
$\gamma_{\mu_1\cdots\mu_k}=\gamma_{[\mu_1}\cdots\gamma_{\mu_k]}$,
e.g.\
$\gamma_{\mu\nu}=\frac{1}{2}\big(\gamma_\mu\gamma_\nu-\gamma_\mu\gamma_\nu\big)$
.

The charge conjugation and the Clif\/ford multiplication give rise
to the well known morphism $\Gamma S\otimes\Gamma
S\hookrightarrow  \Lambda\X(M)$, compare
\cite{Klinker4,Papado2,AlekCortDevProey}. We denote the projection
$\Gamma S\otimes \Gamma S\rightarrow \Lambda^k\X(M)$ by $C_k$ and
its symmetry by $\Delta_k\in\{\pm1\}$. The projection is
explicitly given by the $k$-form
\begin{equation}\label{projectionC}
(C_k(\phi\otimes\psi))_{\mu_1\ldots\mu_k}=C(\phi,\gamma_{\mu_1\cdots
\mu_k}\psi).
\end{equation}
The symmetry of the morphisms obeys $\Delta_k=-\Delta_{k-2}$ and
so may be written as\footnote{This can be made more explicit by
evaluating $\Delta_k$ for $k=0,1$, compare \cite{Klinker4}.}
\begin{equation}\label{delta}
\Delta_k=(-)^{\frac{k(k-1)}{2}} \Delta_0^{k+1} \Delta_1^{k}.
\end{equation}
The charge conjugations as well as the Clif\/ford multiplication
$\gamma:\X(M)\otimes \Gamma S\to \Gamma S$ are parallel with
respect to the Levi-Civita connection  and so are all maps $C_k$.
The map $\Gamma\, {\rm End}(S) \hookrightarrow\Lambda\X(M)$ is
called Fierz relation and an isomorphism onto the image is
explicitly given by
\[
\Omega \mapsto  2^{-[\frac{D}{2}]}  \sum_{n=0}^{\langle D\rangle}
(-)^{\frac{n(n-1)}{2}} \frac{1}{n!}
\,{\text{tr}}(\gamma^{(n)}\Omega) \gamma_{(n)},
\]
with $\langle\dim M\rangle :=\dim M$  if $\dim M$ is even and
$\langle D\rangle :=\frac{1}{2}(\dim M-1)$ if $D$ is odd, compare
\cite{vanProeyen1,Kennedy}. If we take into account the charge
conjugation to identify $S$ and $S^*$ and use (\ref{delta}), the
Fierz identity is written as
\begin{equation}\label{Fierz}
\phi\otimes \psi =\frac{1}{\dim S}\sum\nolimits_{n}
\frac{\Delta_0(\Delta_0\Delta_1)^{n}}{n!} C(\phi,\gamma^{(n)}\psi
) (C\gamma_{(n)}).
\end{equation}

We will often use the notations
\begin{equation*}
\big\{\varphi,\psi\big\}:= 2C_1(\varphi\otimes\psi),\qquad
\langle\phi,\psi\rangle:=C(\phi,\psi).
\end{equation*}
Charge conjugations with $\Delta_1=1$ are of special interest,
because $\big\{\cdot,\cdot\big\}$ may be seen as a supersymmetry
bracket in this case. In particular, this choice is  possible for
Lorentzian space-times, i.e.\ spin manifold of signature
$(-1,1,\ldots,1)$, compare  \cite{Klinker4}. Furthermore, we draw
the attention to \cite{AlekCor2} for a classif\/ication of
bilinear forms also for the case of extended supersymmetry
algebras.
\begin{remark}\label{remarktwist}
Even in the case of $\Delta_1=-1$ we may construct a graded
manifold with  supersymmetry bracket by taking the direct sum of
the spinor bundle with itself and provide it with a~modif\/ied
charge conjugation $C\otimes\tau_2$. Although  there is a choice
of charge conjugation with the appropriate symmetry, we are
sometimes forced to use the ``wrong'' one. For example when we
want to deal with real spinors. We will discuss  such a
construction in section \ref{admissibleconnection} so that we will
omit it here.
\end{remark}

{\em Special vector fields on $\hat M$. } The charge conjugation
$C$ yields an  identif\/ication $S^*\simeq S$. Using this
identif\/ication a natural inclusion $\jmath:\Gamma S\to\X(\hat
M)_1$ is given by the interior multiplication of forms and its
image is a vector f\/ield of degree $-1$. Explicitly we have
\begin{equation*}
 \jmath(\phi): \Gamma S\to \Gamma\Lambda S,\quad \jmath(\phi)(\eta)=\langle\phi,\eta\rangle
\end{equation*}
with the extension as derivation of degree $-1$.

Let us consider a connection $D$ on the spinor bundle $S$. For
every vector f\/ield $X\in\X(M)$ the action of $D_X$ on $\Lambda
S$ is of degree zero. This connection gives rise to an inclusion
$\jmath_D:\X(M)\to\X(\hat M)_0$ given by
\begin{equation}\label{inclusionvectors}
\jmath_D(X):\Gamma\Lambda S\to\Gamma\Lambda S,\qquad
\jmath_D(X)(\eta)=D_X\eta.
\end{equation}

These two inclusions give the  natural splitting
\begin{equation*}
\X(\hat M) = \Gamma\Lambda  S\otimes\X(M)\oplus \Gamma\Lambda
S\otimes\Gamma  S.
\end{equation*}

The endomorphisms of $S$ are vector f\/ields of degree zero on
$\hat M$ in the natural way. Suppose $\Phi \in \Gamma\, {\rm
End}(S) \subset \X(\hat M)_0\cap\X_{\text{v}}(\hat M)$, then the
action is given by
\[
\Gamma\Lambda S\supset \Gamma S\ni \eta \longmapsto  \Phi(\eta)
\in \Gamma S\subset\Gamma \Lambda S.
\]
With respect to a local frame $\{\theta_k\}$ of $S$ the
endomorphism $\Phi$ has the components $\Phi_j^i$ and the
associated vector f\/ield is given by
$\Phi=\Phi_i^jC^{ik}\theta_j\otimes\jmath(\theta_k)$, where
$C^{ij}C_{jk}=\delta^i_k$ and $C_{ij}=C(\theta_i,\theta_j)$.

For  $X,Y \in\X(M)$, $\varphi,\psi\in \Gamma S$ and
$\Phi\in\Gamma \,{\rm End}(\Gamma S)$ the following fundamental
commutation relations hold:
\begin{gather*}
\big[\jmath_D(X),\jmath_D(Y)\big]       = R(X,Y)+\jmath_D([X,Y]), \\
\big[\jmath(\varphi),\jmath(\psi)\big]      = 0,   \qquad
\big[\jmath_D(X),\jmath(\varphi)\big]        = \jmath(D^C_X\phi), \\
\big[\Phi,\jmath(\varphi)\big]      = \jmath(-\Phi^C\varphi),
\qquad \big[\jmath_D(X),\Phi\big]           = D_X\Phi .
\end{gather*}

Consider the space $\X(M)\otimes\Gamma S$ of vector-spinors. The
decomposition into irreducible representation spaces yields
$\X(M)\otimes \Gamma S=\Gamma S\oplus \Gamma S_{\frac{3}{2}}$.
Using  the identif\/ication $\X(M)\simeq\Omega^1(M)$ via~$g$, the
inclusion of $ \Gamma S\hookrightarrow \X(M)\otimes\Gamma S$ is
given  by the Clif\/ford multiplication $\xi(Y)=Y\xi$. In this way
the spin-$\frac{3}{2}$ f\/ields are given by the kernel of the
Clif\/ford multiplication. Given a frame $\{e_\mu\}$ on $M$ with
associated $\gamma$-matrices $\gamma_\mu$ the inclusion is given
by
\[
\Gamma S\hookrightarrow \Gamma S\otimes\X(M),\qquad \phi\mapsto
(\dim M)^{-1}\,\gamma^\mu\phi\otimes e_\mu.
\]
This identif\/ication of the spinors in the vector-spinors is used
to def\/ine a v-like vector f\/ield of degree one on the graded
manifold. For $\phi\in\Gamma S$ we denote this vector f\/ield by
$\imath_D(\phi)$ and it is def\/ined by the above formula up to
the dimension dependent factor together with
(\ref{inclusionvectors}):
\begin{equation*}
\imath_D(\phi)=\gamma^\mu\phi\otimes \jmath_D(e_\mu).
\end{equation*}
In \cite{Klinker2,Klinker4} we used this map with $D=\nabla$ the
Levi-Civita connection on $M$ and $S$. In \cite{Papado2} this
object is considered to construct a (spinor dependent)
dif\/ferential  on $\Lambda S$. The action of the dif\/ferential
corresponds to the action of the vector f\/ield $\imath_D(\phi)$
on the (super)functions $\Gamma\Lambda S$ of~$\hat M$, i.e.
\[
\Gamma\Lambda S\supset\Gamma S \ni \eta  \longmapsto
\imath_(\phi)\eta
    =\gamma^\mu\phi\wedge D_\mu\eta \in \Gamma\Lambda^2S\subset\Gamma\Lambda S.
\]
This vector f\/ield will be considered in section
\ref{JacobiBianchi}.

\section{The torsion of spinor connections}\label{sectionTorsion}

Given a connection $D$ on $S$ we associate to $D$ the f\/ield
$\mathcal{A}:=D-\nabla\in\Omega^1(M)\otimes\Gamma\,{\rm End}(S)$
where as before $\nabla$ denotes the Levi-Civita connection on
$M$. Furthermore, if we denote by~$A$ the projection of
$\mathcal{A}$ onto the sub algebra which is locally given by $
{\rm span}\big\{\gamma_{\mu\nu}\big\}\subset \{
\Phi\in\Gamma\,{\rm End}(S) | [\Phi,\gamma_\mu]\subset {\rm
span}\{\gamma_\nu\} \text{ for all }\mu\} $, then the connection
$\nabla^D=\nabla+A$ is a metric connection on $M$.

As noted in the last section, the charge conjugation $C: S\to S^*
$ as well as the Clif\/ford multiplication $\gamma: TM\otimes
S\to S$ are parallel with respect to the Levi-Civita connection.
More precisely we have the following well known result:

\begin{proposition}
The Clifford multiplication is parallel with respect to the
connection $D$ on $S$ and $\widetilde \nabla $ on $M$ if and only
if $D=\widetilde\nabla$ is a metric connection.

The charge conjugation is parallel with respect to the connection
$D$ on $S$ if and only if $\mathcal{A}$ takes its values in
\[
{\rm span}\big\{\gamma^{\mu_1\cdots \mu_{k}}; \Delta_k\Delta_0=-1
\big\} ={\rm span}\big\{\gamma^{\mu_1\cdots
\mu_{4k+2}},\gamma^{\mu_1\cdots \mu_{4k-\Delta_0\Delta_1}} \big\}
.
\]
In particular, $C$ is parallel with respect to every metric
connection.
\end{proposition}

\begin{example} \label{supergravitation}
{\rm In 11-dimensional space-time, i.e.\ $t=1,s=10$ we have
$\Delta_1=-\Delta_0=1$ so that the  map $\Phi\mapsto\Phi^C, \
\text{with}\  C(\Phi^C\eta,\xi):=C(\eta,\Phi\xi)$ has
($-1$)-eigenspace
\[
{\rm span}\big\{ \gamma^{\mu_1\cdots
\mu_{4k+2}},\gamma^{\mu_1\cdots \mu_{4k+1}}  \big\}
\]
and ($+1$)-eigenspace
\[
{\rm span}\big\{ \gamma^{\mu_1\cdots
\mu_{4k+3}},\gamma^{\mu_1\cdots \mu_{4k}}\big\}.
\]
In particular, the Clif\/ford multiplication is skew symmetric.
For example, consider  the supercovariant derivation which come
from the supergravity variation of the gravitino and for which
$\mathcal{A}_X$ has a three-form and a f\/ive-form part, compare
\cite{CremmJulScher,CremmJul,PapadoTsi,GreenSchwarzWitten}. This
connection does not make the charge conjugation parallel. }
\end{example}
Due to this example, parallelism of the charge conjugation is not
the appropriate notion to be related to supersymmetry in general.

To the connection  $D$ on $S$ we will associate  another
connection $D^C$. To construct this  we consider the connection
$D\otimes \mathbbm{1}+\mathbbm{1}\otimes D^C$  on $S\otimes S$ and
the induced connection on  $S^*\otimes S^*$. Then $D\otimes
\mathbbm{1}+\mathbbm{1}\otimes D^C$ shall make the charge
conjugation parallel, i.e.\ $(D\otimes
\mathbbm{1}+\mathbbm{1}\otimes D^C)C=0$. For
$D=\nabla+\mathcal{A}$ this implies $D^C=\nabla-\mathcal{A}^C$.
The next remark is obtained immediately.
\begin{remark}\label{curvC}
The curvature $R$ of the connection $D$ and the curvature $R^C$ of
the connection $D^C$ are related by $(R(X,Y))^C=-R^C(X,Y)$.
\end{remark}
We endow the bundle of ${\rm End}(S)$-valued tensors on $M$ with a
connection induced by $D$, $D^C$ and $\nabla$.
\begin{definition}\label{definitionhatD}
Let $\Phi\in \X(M)^{\otimes k} \otimes \Omega^1(M)^{\otimes\ell}
\otimes \Gamma\, {\rm End}(S)$. The connection $\hat D$ is
def\/ined by
\begin{equation*}
(\hat D_Z\Phi)(X)\xi := D_Z(\Phi(X)\xi)-\Phi(\nabla
_ZX)\xi-\Phi(X)D^C_Z\xi
\end{equation*}
for all  vector f\/ields $Z, X\in \Omega^1(M)^{\otimes k}\otimes
\X(M)^{\otimes\ell}$, and  $\xi\in\Gamma S$.
\end{definition}
We consider the following ${\rm ad}$-type representation of ${\rm
End}(S)$ on itself.
\begin{definition}
Let $\Omega\in{\rm End}(S)$. We def\/ine ${\rm ad}^C_\Omega: {\rm
End}(S)\to {\rm End}(S)$ by
\begin{equation*}
{\rm ad}^C_\Omega \Phi:= \Omega\Phi+\Phi \Omega^C.
\end{equation*}
\end{definition}
This is indeed a representation, because  ${\rm
ad}^C_{[\Omega_1,\Omega_2]}
 \Phi= \big[{\rm ad}^C_{\Omega_1},{\rm ad}^C_{\Omega_2}\big]\Phi$.
For $\Omega=\Omega^++\Omega^-$, i.e.\ $\Omega^C=\Omega^+-\Omega^-$
we have
\begin{equation*}
{\rm ad}^C_\Omega \Phi=
\big[\Omega_-,\Phi\big]+\big\{\Omega_+,\Phi\big\} .
\end{equation*}
Furthermore we have
\begin{equation}\label{pres}
({\rm ad}^C_\Omega \Phi)^C =
(\Omega\Phi+\Phi\Omega^C)^C=\Phi^C\Omega^C+\Omega\Phi^C={\rm
ad}^C_\Omega\Phi^C
\end{equation}
which yields
\begin{proposition}\label{eigenpreserve}
${\rm ad}^C_\Omega$ preserves the $(\pm 1)$-eigenspaces of the
linear map $\Phi\mapsto\Phi^C$ for all $\Omega\in{\rm End}(S)$.
\end{proposition}

\begin{proposition}
Let $\hat D$ be the connection associated to $D$ cf.\ Definition
{\rm \ref{definitionhatD}}. Then $D$ and the charge adjoint are
compatible in the way that
\begin{equation}\label{compatibleDad}
\hat D({\rm ad}^C_\Omega\Psi)={\rm ad}^C_{D\Omega}\Psi+{\rm
ad}^C_\Omega\hat D\Psi
\end{equation}
for all $\Omega,\Psi\in \Gamma \,{\rm End}(S)$.
\end{proposition}

\begin{definition}\label{ctorsion}
Let $D$ be a connection on the spinor bundle $S$ over the (pseudo)
Riemannian manifold $M$ and denote by $\nabla$  the Levi-Civita
connection on $M$. The {\sc torsion}
$\mathcal{T}\in\Omega^2(M)\otimes\Gamma\,{\rm End}(S)$ of  $D$ is
the def\/ined by two times the skew symmetrization of $\hat
D\gamma:\X(M)\otimes\X(M)\to \Gamma \,{\rm End}(S)$.
\end{definition}

\begin{remark} \ \ {}

\begin{enumerate}\itemsep=0pt
\item We have $(\hat D_X\gamma)(Y)=\hat
D_X(\gamma(Y))-\gamma(\nabla_XY)$.
 Using this and $\nabla_XY-\nabla_YX=[X,Y]$ and omitting  the map $\gamma$ we may also write
\[
\mathcal{T}(X,Y)=\hat D_XY-\hat D_YX-[X,Y] .
\]
\item In terms of the dif\/ference $\mathcal{A}=D-\nabla
\in\Omega^1(M)\otimes\Gamma \, {\rm End}(S)$ the torsion may be
written as
\[
\mathcal{T}(X,Y)= {\rm ad}^C_{\mathcal{A}(X)}Y-{\rm
ad}^C_{\mathcal{A}(Y)}X.
\]
\item The last point and (\ref{pres}) yield that the torsion has
symmetry $\Delta_1$, i.e.\ for all $\eta$, $\xi$ we have
\[
C(\eta,\mathcal{T}_{\mu\nu}\xi)=\Delta_1C(\xi,\mathcal{T}_{\mu\nu}\eta)\,.
\]
\item For a metric connection $D$ on $S$ the torsion is exactly
the torsion which is def\/ined by the connection $D$ on the
manifold $M$.
\end{enumerate}
\end{remark}

The torsion obeys some  Bianchi-type identities.
\begin{proposition}\label{Bianchi}
Let $D$ be a connection on the spinor bundle $S$ over the
(pseudo) Riemannian manifold $M$. The torsion  $\mathcal{T}$ and
the curvature $R$ of $D$ obey
\begin{gather}
\hat D_{[\kappa}\mathcal{T}_{\mu\nu]}= {\rm ad}^C( R_{[\kappa\mu})\gamma_{\nu]},\label{BianchiDT}\\
\hat D_{[\kappa}({\rm ad}^C_ R\gamma)_{\mu\nu\rho]}={\rm
ad}^C(R_{[\kappa\mu})\mathcal{T}_{\nu\rho]}.\label{BianchiDadR}
\end{gather}
In this context we add the following identity for the curvature
$R$ of $D$\footnote{This identity  holds for any connection $D$ on
a vector bundle over $M$, if we endow all tensor bundles with the
connection induced by $D$ and the Levi-Civita connection on $M$.}:
\begin{equation}
D_{[\kappa}R_{\mu\nu]}=0.\label{BianchiDR}
\end{equation}
\end{proposition}
\begin{proof}
With Def\/inition \ref{definitionhatD} the left hand side of
(\ref{BianchiDT}) is given by
\begin{equation*}
(\hat D_X\mathcal{T})(Y,Z)\xi
 =  D_X(\mathcal{T}(X,Y)\xi)-\mathcal{T}(\nabla_XY,Z)\xi-\mathcal{T}(Y,\nabla_XZ)\xi-\mathcal{T}(Y,Z)D^C_X\xi.
\end{equation*}
We use the def\/inition of the torsion and get
\begin{align*}
(\hat D_X\mathcal{T})(Y,Z)\xi
 =&\ D_X(\mathcal{T}(Y,Z)\xi)-\mathcal{T}(\nabla_XY,Z)\xi- \mathcal{T}(Y,\nabla_XZ)\xi -\mathcal{T}(Y,Z)D^C_X\xi \\
 =&\ D_XD_Y(Z\xi)-\underline{ D_X(\nabla_YZ\xi)}-\underline{D_X(ZD^C_Y\xi)} -D_XD_Z(Y\xi)\\
   & +\underline{D_X(\nabla_ZY\xi)} +\underline{D_X(YD^C_Z\xi)}-D_{\nabla_XY}(Z\xi)+\nabla_{\nabla_XY}Z\xi \\
   & +ZD^C_{\nabla_XY}\xi +\underline{D_Z(\nabla_XY\xi)} -\nabla_Z\nabla_XY\xi-\underline{\nabla_XYD^C_Z\xi} \\
   & -\underline{D_Y(\nabla_XZ\xi)}+\nabla_Y\nabla_XZ\xi+\underline{\nabla_XZD^C_Y\xi}+D_{\nabla_XZ}(Y\xi) \\
   &  -\nabla_{\nabla_XZ}Y\xi -YD^C_{\nabla_XZ}\xi -\underline{D_Y(ZD^C_X\xi)}+\underline{\nabla_YZD^C_X\xi} \\
   & +ZD^C_YD^C_X\xi + \underline{D_Z(YD^C_X\xi)}-\underline{\nabla_ZYD^C_X\xi}-YD^C_ZD^C_X\xi .
\end{align*}
The underlined terms are symmetric with respect to $X$, $Z$ or
$X$, $Y$. So they vanish when we skew symmetrize the above
expression with respect to $X$, $Y$, $Z$. So we are left with
\begin{gather*}
  (\hat D_X\mathcal{T})(Y,Z)\xi   + (\hat D_Z\mathcal{T})(X,Y)\xi+(\hat D_Y\mathcal{T})(Z,X)\xi  \\
  \qquad{}  =    R(X,Y)(Z\xi)+R(Z,X)(Y\xi)+R(Y,Z)(X\xi)  +Z R^C(Y,X)\xi + YR^C(X,Z)\xi \\
  \qquad\phantom{=}{} + XR^C(Y,Z)\xi + \underbrace{(R^0(Y,X)Z+R^0(Z,Y)X+R^0(X,Z)Y )}_{=0}\xi.
\end{gather*}
With Remark \ref{curvC}, i.e.\ $R^C(Y,X)=R(X,Y)^C$, we may rewrite
this as
\begin{equation*}
\hat D_{[\mu}\mathcal{T}_{\nu\kappa]} =
R_{[\mu\nu}\gamma_{\kappa]}+\gamma_{[\kappa}(R_{\mu\nu]})^C
                                                 =  {\rm ad}^C(R_{[\mu\nu})\gamma_{\kappa]}.
\end{equation*}
The proof  of (\ref{BianchiDR})  is done by similar calculations.
(\ref{BianchiDadR}) follows from (\ref{compatibleDad}) and
(\ref{BianchiDR}) after skew symmetrization of
\begin{align*}
(\hat D_\kappa {\rm ad}^C_R\gamma)_{\mu\nu\rho}={\rm
ad}^C_{D_\kappa R_{\mu\nu}} \gamma_\rho + {\rm
ad}^C_{R_{\mu\nu}}\hat D_\kappa\gamma_\rho.
\end{align*}
This completes the proof.
\end{proof}

\begin{example}
{\rm We consider a manifold  which admits  geometric Killing
spinors. These are spinors which fulf\/ill the equation
$\nabla_X\phi=-aX\phi$ for a constant $a\neq 0$, the Killing
number. This equation has been extensively examined in the
literature \cite{Baum5,Moroi1,KathHabil} and in particular
\cite{BFGK}. Moreover we would  like to stress on \cite{Baer1}
where the author draws a remarkable connection between geometric
Killing spinors on a manifold and parallel spinors on the cone
over the manifold, at least in the Riemannian case.

From the above equation we read that the connection $D$ on the
spinor bundle for which the geometric Killing spinors are parallel
is given by
\begin{equation*}
D=\nabla+a\cdot\gamma.
\end{equation*}
Suppose $\Delta_1\Delta_0=-1$, i.e.\ the Clif\/ford multiplication
is skew symmetric. This yields a condition on the connection which
will be important in the next section:
\begin{align*}
{\rm ad}^C(\mathcal{A}_{\{\mu})\gamma_{\nu\}}
    & =a\gamma_{\{\mu}\gamma_{\nu\}}+a\gamma_{\{\nu}\gamma^C_{\mu\}}
       =-ag_{\mu\nu}+\Delta_0\Delta_1a\gamma_{\{\nu}\gamma_{\mu\}}\\
    & =-a(1+\Delta_1\Delta_0)g_{\mu\nu}
       =0 .
\end{align*}
The torsion and the curvature of this connection are given by
\[
\mathcal{T}_{\mu\nu}=4a \gamma_{\mu\nu} \qquad \text{and}\qquad
R_{\mu\nu}=R^0_{\mu\nu}+2a^2\gamma_{\mu\nu}
\]
and obey
\[
\hat D_\kappa\mathcal{T}_{\mu\nu}=-16ag_{\kappa[\mu}\gamma_{\nu]}
\quad \text{and}\quad
ad^C_{R_{\mu\nu}}\gamma_\kappa=R^0_{\mu\nu\kappa\lambda}\gamma^\lambda+8a^2g_{\kappa[\mu}\gamma_{\nu]}.
\]
such that both sides of (\ref{BianchiDT}) vanish. }
\end{example}

\section{Admissible spinor connections}\label{admissibleconnection}

\subsection{Killing equations and admissible connections}

We examine the conditions on the connection
$D=\nabla+\mathcal{A}$, such that the vector f\/ield
$\{\varphi,\psi\}$ built up by the Killing spinors
$D^C\varphi=D^C\psi=0$ is a Killing vector f\/ield, i.e.\
$\mathcal{L}_{\{\varphi,\psi\}}g=0$. We have
\begin{align*}
\mathcal{L}_{\{\varphi,\psi\}}g(e_\mu,e_\nu)
&= g(\nabla_\mu\{\varphi,\psi\},e_\nu)+g(\nabla_\nu\{\varphi,\psi\},e_\mu) \\
&= g(\{\nabla_\mu\varphi,\psi\},e_\nu)+g(\{\varphi,\nabla_\mu\psi\},e_\nu) +\{\mu\leftrightarrow\nu\}\\
&=g(\{\mathcal{A}^C_\mu\varphi,\psi\},e_\nu)+g(\{\varphi,\mathcal{A}^C_\mu\psi\},e_\nu)+\{\mu\leftrightarrow\nu\}\\
&=2\langle  \mathcal{A}^C_\mu\varphi,\gamma_\nu\psi \rangle  +2\langle \varphi,\gamma_\nu\mathcal{A}^C_\mu\psi\rangle+\{\mu\leftrightarrow\nu\}\\
&= 2\langle \varphi, ad^C_{\mathcal{A}_\mu}\gamma_\nu \psi\rangle
+\{\mu\leftrightarrow\nu\}.
\end{align*}
This yields
\begin{theorem}\label{KillingA}
Let $D$ be a connection on the spinor bundle $S$ over $M$. Suppose
$\phi,\psi\in S$ are parallel with respect to the associated
connection $D^C$. Then the vector field $\{\phi,\psi\}=2
C_1(\phi\otimes \psi)$ is a Killing vector field if the symmetric
part of $\hat D\gamma: \X(M)\otimes \X(M) \to {\rm End}(S)$ acts
trivially on the parallel spinors. In this case we have
\[
\nabla_\mu\{\eta,\xi\}_\nu=C(\eta,\mathcal{T}_{\mu\nu}\xi).
\]
\end{theorem}
This motivates the next def\/inition.
\begin{definition}\label{defadmissible}
Let $D$ be a connection on the spinor bundle $S$ over $M$  and
$\mathcal{K}\subset \Gamma S$ be a subset.
\begin{enumerate}\itemsep=0pt
\item We call $(\mathcal{K},D)$ admissible if the symmetric part
of $\hat D\gamma$ acts trivially on $\mathcal{K}$. If $D$ is
f\/ixed we call $\mathcal{K}$ admissible. \item We call $D$
admissible if  $\hat D\gamma$  is skew symmetric. In this case is
$\mathcal{T}=2\hat D\gamma$.
\end{enumerate}
\end{definition}
\begin{remark}
Due to Theorem \ref{KillingA} the admissible subsets of
$D^C$-parallel spinors are of particular interest.
\end{remark}
\begin{example}\label{sugra}
{\rm Consider the supergravity connection $D=\nabla+\mathcal{A}$
with $\mathcal{A}=F^3+F^5$ given by
\[
F^3(X)=-\frac{1}{36}X^\mu
F_{\mu\nu\rho\sigma}\gamma^{\nu\rho\sigma} \qquad \text{and}\qquad
F^5(X)=\frac{1}{288}X_\mu
F_{\nu\rho\sigma\tau}\gamma^{\mu\nu\rho\sigma\tau}
\]
for a 4-form $F$ on $M$. This connection obeys
$F_\mu^5=-(F_\mu^5)^C$ and $F_\mu^3=(F_\mu^3)^C$ due to Example
\ref{supergravitation}. Furthermore we have
\begin{align*}
{\rm ad}^C(\mathcal{A}_\mu)\gamma_\nu & =
\big[F^5_\mu,\gamma_\nu\big]+\big\{F^3_\mu,\gamma_\nu\big\}
  = \frac{1}{144} F^{\kappa\rho\sigma\tau} \gamma_{\mu\nu\kappa\rho\sigma\tau}
     +\frac{1}{9} F_{\mu\nu\kappa\rho} \gamma^{\kappa\rho}
\end{align*}
which is indeed skew symmetric with respect to $\mu$ and $\nu$,
i.e.\ the supergravity connection is admissible. }
\end{example}
This example can be generalized.

\begin{theorem}\label{admissibleForm}
Let $D$ be a connection on the spinor bundle $S$ of $M$ and
$\mathcal{A}:=D-\nabla\in\Omega^1(M)\otimes \Gamma\, {\rm
End}(S)$. Suppose $\mathcal{A}_X$ is homogeneous  with respect to
$\Gamma\,{\rm End}(S) \simeq\bigoplus_k\Omega^k(M)$. Consider the
decomposition\footnote{$\Omega^{(\ell,1)}$ denotes the irreducible
representation space with highest weight $e_1+e_\ell$.}
\[
\Omega^1(M)\otimes\Omega^\ell(M)
=\Omega^{\ell+1}(M)\oplus\Omega^{\ell-1}(M)\oplus\Omega^{(\ell,1)}.
\]
$\mathcal{A}(X)$ may be written as $\mathcal{A}_X= X \rfloor
F^{\ell+1}+X\wedge G^{\ell-1}+\mathcal{A}_0(X)$ with an
$(\ell+1)$-form $F$, an $(\ell-1)$-form $G$ and
$\mathcal{A}_0\in\Omega^{(\ell,1)}$.

Then $D$ is admissible if and only if $\mathcal{A}_0=0$ and
$\Delta_1\Delta_{\deg}=-1$ or equivalently
$\Delta_0\Delta_{{\deg}-1}=(-)^{\deg}$, i.e.\ ${\deg}\equiv
3\,{\text{\rm mod}}\, 4$, or $1+\Delta_0\Delta_1\,{\text{\rm
mod}}\, 4$. Here ${\deg}$ denotes the degree of the forms~$F$
and~$G$ respectively.
\end{theorem}

\begin{proof}
Consider $\mathcal{A}$ either to be of the form
\[ A_\mu=F_{\mu\kappa_1\ldots\kappa_\ell}\gamma^{\kappa_1\ldots\kappa_\ell}
\quad\text{or}\quad
A_\mu=G^{\kappa_1\ldots\kappa_{\ell-1}}\gamma_{\mu\kappa_1\ldots\kappa_{\ell-1}}
\]
with $F\in\Omega^{\ell+1}\oplus\Omega^{(\ell,1)}$ and
$G\in\Omega^{\ell-1}$. In the f\/irst case we have
\begin{gather*}
\mathcal{T}_{\mu\nu}
    = F_{\mu}{}^{\kappa_1\ldots\kappa_\ell} \big( \gamma_{\kappa_1\ldots\kappa_\ell}\gamma_\nu
        +\gamma_\nu \gamma_{\kappa_1\ldots\kappa_\ell}^C\big) \\
 \phantom{\mathcal{T}_{\mu\nu}}{}   = F_{\mu}{}^{\kappa_1\ldots\kappa_\ell} \big( \gamma_{\kappa_1\ldots\kappa_\ell}\gamma_\nu
        +\Delta_0\Delta_\ell\gamma_\nu \gamma_{\kappa_1\ldots\kappa_\ell}\big) \\
 \phantom{\mathcal{T}_{\mu\nu}}{}    = F_{\mu}{}^{\kappa_1\ldots\kappa_\ell} \big( \gamma_{\kappa_1\ldots\kappa_\ell\nu}
        +\Delta_0\Delta_\ell \gamma_{\nu\kappa_1\ldots\kappa_\ell}\big) 
        -\ell F_{\mu}{}^{\kappa_1\ldots\kappa_\ell} \big((-)^{\ell-1}+\Delta_0\Delta_\ell \big)
        g_{\nu[\kappa_1} \gamma_{\kappa_2\ldots\kappa_\ell]} \\
 \phantom{\mathcal{T}_{\mu\nu}}{}    = F_{\mu}{}^{\kappa_1\ldots\kappa_\ell} \big( 1+(-)^\ell\Delta_0\Delta_\ell\big) \gamma_{\kappa_1\ldots\kappa_\ell\nu}
        -\ell F_{\mu\nu}{}^{\kappa_2\ldots\kappa_\ell} \big((-)^{\ell-1}+\Delta_0\Delta_\ell \big)
         \gamma_{\kappa_2\ldots\kappa_\ell}.
\end{gather*}
This expression is skew symmetric if and only if $F$ is totally
skew symmetric, i.e.\ $\mathcal{A}_0=0$, and
$\Delta_0\Delta_\ell=(-)^{\ell-1}$. With ${\deg}=\ell+1$ this is
exactly the condition stated. The second case is treated in almost
the same way.
\begin{gather*}
\mathcal{T}_{\mu\nu}
    = G^{\kappa_1\ldots\kappa_{\ell-1}} \big( \gamma_{\mu\kappa_1\ldots\kappa_{\ell-1}}\gamma_\nu
        +\gamma_\nu \gamma_{\mu\kappa_1\ldots\kappa_{\ell-1}}^C\big) \\
\phantom{\mathcal{T}_{\mu\nu}}{}    =
G^{\kappa_1\ldots\kappa_{\ell-1}} \big(
\gamma_{\mu\kappa_1\ldots\kappa_{\ell-1}}\gamma_\nu
        +\Delta_0\Delta_\ell\gamma_\nu \gamma_{\mu\kappa_1\ldots\kappa_{\ell-1}}\big) \\
\phantom{\mathcal{T}_{\mu\nu}}{}     =
G^{\kappa_1\ldots\kappa_{\ell-1}} \big(
\gamma_{\mu\kappa_1\ldots\kappa_{\ell-1}\nu}
        +\Delta_0\Delta_\ell \gamma_{\nu\mu\kappa_1\ldots\kappa_{\ell-1}}\big) 
        -\ell G^{\kappa_1\ldots\kappa_{\ell-1}} \big( (-)^{\ell-1}+\Delta_0\Delta_\ell\big)
        g_{\nu[\mu} \gamma_{\kappa_1\ldots\kappa_{\ell-1}]}\! \\
\phantom{\mathcal{T}_{\mu\nu}}{}    =
G^{\kappa_1\ldots\kappa_{\ell-1}} \big(
(-)^{\ell-1}+\Delta_0\Delta_\ell\big)
        \gamma_{\mu\nu\kappa_1\ldots\kappa_\ell}  
        -\ell G^{\kappa_1\ldots\kappa_{\ell-1}} \big( (-)^{\ell-1}+\Delta_0\Delta_\ell\big)
        g_{\nu[\mu} \gamma_{\kappa_1\ldots\kappa_{\ell-1}]}.
\end{gather*}
This is skew symmetric if and only if
$\Delta_0\Delta_\ell=(-)^\ell$ or
$\Delta_0\Delta_{\ell-2}=(-)^{\ell-1}$ which with ${\deg}=\ell-1$
f\/inishes the proof.
\end{proof}

If $\mathcal{A}$ is of the form $\mathcal{A}_X =\hat\alpha X\wedge
F +\hat\beta X\rfloor F$ for an $\ell$-form $F$ we may rewrite it
as $\mathcal{A}_X=\alpha X\cdot F+\beta F\cdot X$ where $(\cdot)$
denotes Clif\/ford multiplication and $\alpha$, $\beta$ are linear
combinations of~$\hat\alpha$,~$\hat\beta$. Therefore, we will
restrict ourself often to the two cases $F\cdot X$ and $X\cdot F$.

\begin{remark} \ \ \null

\begin{itemize}\itemsep=0pt
\item To be admissible is a property which has to be checked for
every degree of ${\rm ad}^C_{\mathcal{A}}\gamma$. This yields that
the connection $D$ on $S$ is admissible if and only if every
homogeneous summand is. Furthermore $D$ is admissible if\/f $D^C$
is admissible, because this fact does only depend on the degree of
$\mathcal{A}_\mu$ in $\Omega(M)$ which  is independent of the
charge conjugation. \item For $\mathcal{A}_X=X\wedge F +X\rfloor
G$ admissible the torsion  is given by
\begin{equation*}
\mathcal{T}(X,Y)= \pm X\wedge Y\wedge F \pm X\rfloor Y\rfloor G .
\end{equation*}
\end{itemize}
\end{remark}

\begin{example}
{\rm Let $\mathcal{A}(X)$ be of the form $X\wedge F^{(\ell)}$ or
$X\rfloor F^{(\ell)}$. In  eleven dimensional space time this
leads to an admissible connection for $\ell=0,3,4,7,8,11$. In
Example \ref{sugra} we have $F^{3}_X\sim X\rfloor F^{(4)}$ and
$F^{5}_X\sim X\wedge F^{(4)}$. }
\end{example}

Theorems \ref{KillingA} and \ref{admissibleForm} have an important
consequence for metric connections on the spinor bundle.
\begin{corollary}
Let $D$ be a metric connection on $S$. $D$ is admissible if and
only if $\mathcal{A}_X$  is of the form $X\rfloor F^{(3)}$. We
write  $\mathcal{A}_\mu
=\frac{1}{4}A_{\mu\nu\kappa}\gamma^{\nu\kappa}$. The torsion
tensor in this case is totally skew symmetric and given by
$T_{\mu\nu\kappa}=2A_{[\mu\nu]\kappa}= 2A_{\mu\nu\kappa}$. In
other words $D$ is admissible if and only if its torsion is
totally skew symmetric.
\end{corollary}
Metric connections with skew symmetric torsion play an important
role in string theory as well as supergravity theories. A lot of
literature on this topic has been published during the past few
years, see for example \cite{Ivanov1} or~\cite{IvanovIvanov} and
references therein.

\subsection{Admissible connections on twisted spinor bundles}

Sometimes it is necessary to introduce $\ell$-form f\/ields which
have  degree dif\/ferent from those which are allowed by
Theorem~\ref{admissibleForm}. This is possible in two dif\/ferent
ways.

The f\/irst way is, in particular, interesting if $M$ is  of even
dimension $2n$.

Suppose $n$ is even. In this case the $\ell$-forms with
$\ell\equiv 1$ or $1+\Delta_0\Delta_1\,{\text{mod}}\,4$ contribute
to an admissible connection by
\begin{equation}\label{starcontributioneven}
F^{\nu_1\ldots\nu_\ell}\gamma_{\nu_1\ldots\nu_\ell}\gamma_\mu
\gamma^*.
\end{equation}
This is due to $\Delta(\gamma^{(\ell)}\gamma^*)=\Delta_{2n-\ell}$
(compare (\ref{starsymmetry}) in Appendix \ref{gammaappendix}) and
\[
2n-1\equiv 3\,{\text{mod}}\,4,\qquad 2n-3\equiv
1\,{\text{mod}}\,4,\qquad 2n-(1\pm\Delta_1\Delta_0)\equiv 1\pm
\Delta_1\Delta_0\,{\text{mod}}\,4
\]
for $n$ even as well as Theorem~\ref{admissibleForm}.

If $n$ is odd  we have
\[
2n-1\equiv 1\,{\text{mod}}\,4,\qquad 2n-3\equiv
3\,{\text{mod}}\,4,\qquad 2n-(1\pm\Delta_1\Delta_0)\equiv 1\mp
\Delta_1\Delta_0\,{\text{mod}}\,4.
\]
In this case the $\ell$-form with $\ell\equiv 3$ or
$1-\Delta_0\Delta_1 \,{\text{mod}}\,4$ contributes cf.\
(\ref{starcontributioneven}), for the same reason.

The introduction of $\gamma^*$ is a bit artif\/icial, because we
may express for example $
F^{\nu_1\ldots\nu_\ell}\gamma_{\nu_1\ldots\nu_\ell}\cdot X
\gamma^*$ as $\pm(*
F)^{\nu_1\ldots\nu_{2n-\ell}}\gamma_{\nu_1\ldots\nu_{2n-\ell}} X$.
Nevertheless, we will see  in Section~\ref{examplePure} that this
is a useful description.
\begin{corollary}\label{restrictive}
We consider the projections $\Pi^\pm:S=S^+\oplus S^-\to S^\pm$. An
$\ell$-form contributes to an admissible connection by
$F_{(\ell)}\gamma^{(\ell)}\gamma_\mu\Pi^\pm$ if and only if
$\ell\equiv 3\,{\text{\rm mod}}\,4$ for $n$ odd, or $\ell\equiv
1+\Delta_0\Delta_1\,{\text{\rm mod}}\,4$ for $n$ even.
\end{corollary}

The second way uses the forms without considering duality, i.e.\
without adding $\gamma^*$. This bypasses the last remark.

We replace the spinor bundle $S$ by the direct sum $S\oplus S$.
This space is equipped with a charge conjugation which is given by
the charge conjugation $C$ on  $S$  twisted by a modif\/ied
Pauli-matrix, i.e.\ $C\otimes \tau_i$. For $\tau_0$ we get the
direct sum of $C$ and we denote this usually by $C$, too. The
connection $D$ for an $\ell$-form $F$ may  be written as
\[
D_\mu=\nabla_\mu+\tfrac{1}{\ell!}F_{i_1\ldots
i_\ell}\gamma^{i_1\ldots
i_\ell}\gamma_\mu\otimes\tau_j=\nabla_\mu+F \gamma_\mu\tau_j\,
\]
with a matrix $\tau_j$. We have
\begin{gather*}
C\otimes \tau_i(F\gamma_\mu\tau_j\eta ,\gamma_\nu\xi)  =
C(F\gamma_\mu\tau_j\eta ,\gamma_\nu\tau_i\xi)
 =  \Delta_\ell\Delta_1\varepsilon_jC(\eta, \gamma_\mu F\gamma_\nu\tau_j\tau_i\xi) \\
\phantom{C\otimes \tau_i(F\gamma_\mu\tau_j\eta ,\gamma_\nu\xi)}{}
 =  \Delta_\ell\Delta_1\varepsilon_j\varepsilon_{ij}C\otimes \tau_i (\eta, \gamma_\mu F\gamma_\nu\tau_j\xi).
\end{gather*}
This yields
\begin{theorem}\label{admissibleForm2}
For the twisted spinor bundle $S\oplus S$ with charge conjugation
$C\otimes\tau_i$ the $\ell$-form $F$ contributes to an admissible
connection in the form $\tfrac{1}{\ell!}F_{i_1\ldots
i_\ell}\gamma^{i_1\ldots i_\ell}\gamma_\mu\otimes\tau_j$ if and
only if
\begin{equation*}
\Delta_\ell\Delta_1\varepsilon_j \varepsilon_{ij}=-1.
\end{equation*}
All possible values for $(\ell,i,j)$ are listed in Table~{\rm
\ref{tableList}}.
\end{theorem}

If we f\/ix $\ell$ we see that the possible values of $j$  depend
on the choice of
 $\tau_i$ in the charge conjugation and on $\Delta_0\Delta_1$ (at least for even $\ell$).
For $i=j=0$ the two components decouple and we recover the result
from Theorem \ref{admissibleForm}.

\begin{table}[htb]\caption{Possible choices for $\tau_j$ so that the $\ell$-form
contributes to an admissible connection if the charge conjugation
is given by $C\otimes \tau_i$.}\label{tableList}
$$
\begin{array}{c|c|c}
i   & \ell\,{\text{mod}}\,4     & j     \\\hline\hline 0   & 1
& 2 \\\hline
    & 3         & 0,1,3 \\\hline\hline
1   & 1         & 3 \\\hline
    & 3         & 0,1,2 \\\hline\hline
2   & 1         & 1,2,3 \\\hline
    & 3         & 0 \\\hline\hline
3   & 1         & 1 \\\hline
    & 3         & 0,2,3
\end{array}
 \qquad \qquad
\begin{array}{c|c|c}
i   & \ell\,{\text{mod}}\,4     & j     \\\hline\hline 0   &
1-\Delta_0\Delta_1    & 2 \\\hline
    & 1+\Delta_0\Delta_1    & 0,1,3 \\\hline\hline
1   & 1-\Delta_0\Delta_1    & 3 \\\hline
    & 1+\Delta_0\Delta_1    & 0,1,2 \\\hline\hline
2   & 1-\Delta_0\Delta_1    & 1,2,3 \\\hline
    & 1+\Delta_0\Delta_1    & 0 \\\hline\hline
3   & 1-\Delta_0\Delta_1    & 1 \\\hline
    & 1+\Delta_0\Delta_1    & 0,2,3
\end{array}
$$
\end{table}

\begin{remark}
We draw the attention to the fact that we change the symmetry of
$C_1$ if we use $\tau_2$ to modify the charge conjugation.
\end{remark}
\begin{example}\label{exampleBergRoo}
\ \ \null

\begin{itemize}\itemsep=0pt
\item In \cite{AliGanGhoPar} the authors discuss pp-wave solutions
of  type IIA supergravity. The starting point is a Killing
equation for the spinors constructed by a 3-form  and a 4-form  in
the following way
\[
D= \nabla+F^3\gamma\otimes\tau_3+F^4\gamma\otimes\tau_1.
\]
In ten dimensional space time we have two natural ways to choose
the charge conjugation ($\Delta_0=+1$ or $-1$) and in both cases
we have $\Delta_1=1$. The above connection is admissible for the
choice $\Delta_0=-1$ and $i=0$ or $3$ as we read from
Table~\ref{tableList}. \item In \cite{BergRooJanOrt} the type IIB
supergravity and the variation of its f\/ields are discussed. The
vanishing of the gravitino variation leads to a Killing equation
where $\mathcal{A}^C$ contains all odd $\ell$-forms, $F^\ell$
which are twisted by $\tau_2$ if the degree is $\ell\equiv
1\,{\text{mod}}\, 4$ and by $\tau_1$ if the degree is $\ell\equiv
3\,{\text{mod}}\, 4$ and furthermore a second three-form, $H^3$,
twisted by $\tau_3$. This is possible only for $i=0$ independent
of $\Delta_0$.

Moreover the four $\z_2$-symmetries  which are given by
multiplying the fermion doublets by~$\tau_1$ or $\tau_3$  may be
seen as a change of the charge conjugation from $C\otimes\tau_0$
to $C\otimes\tau_1$ or $C\otimes\tau_3$. Now it is evident from
Theorem~\ref{admissibleForm2} that not all of\/f f\/ields are
allowed if we want to keep the connection admissible. In
particular, these are $F^\ell=0$ for all $\ell$ if $j=1$ and
$F^1=F^5=F^9=H^3=0$ if $j=3$. These are exactly the truncations
which are made in~\cite{BergRooJanOrt}.
\end{itemize}
\end{example}

We carry on considering  the supergravity connection cf.\
\cite{BergRooJanOrt} which is  given by $D^C=\nabla-\mathcal{A}^C$
with
\begin{gather}
\mathcal{A}^C_\mu = H_{\mu\kappa\lambda}\gamma^{\kappa\lambda} +
F^1_{\kappa}\gamma^\kappa\gamma_\mu\otimes\tau_{2}
+\tfrac{1}{3!}F^3_{\kappa_1\kappa_2\kappa_3}\gamma^{\kappa_1\kappa_2\kappa_3}\gamma_\mu\otimes\tau_{1}
 +\tfrac{1}{5!}F^5_{\kappa_1\ldots\kappa_5}\gamma^{\kappa_1\ldots\kappa_5}\gamma_\mu\otimes\tau_{2}\nonumber\\
 \phantom{\mathcal{A}^C_\mu =}{}
  +\tfrac{1}{7!}F^7_{\kappa_1\ldots\kappa_7}\gamma^{\kappa_1\ldots\kappa_7}\gamma_\mu\otimes\tau_{1}
+\tfrac{1}{9!}F^9_{\kappa_1\ldots\kappa_9}\gamma^{\kappa_1\ldots\kappa_9}\gamma_\mu\otimes\tau_{2},\label{bsp}
\end{gather}
where $H$ is a torsion three form and the $\ell$-forms are
connected by $*F^1=F^9$, $*F^3=-F^7$, and $*F^5=F^5$. As we
mentioned in Example~\ref{exampleBergRoo} this connection is
admissible for the charge conjugation~$C\otimes\tau_0$.

Due to the nature of the gravity theories the parallel spinors
have a f\/ixed chirality property. More precisely the chirality of
the two components of $\eta$ and a relation between both entries,
are f\/ixed for all supersymmetry parameters. This may be
described by an operator
\begin{equation}\label{projection}
\mathbbm{1}\otimes\tau_i\pm\gamma^*\otimes\tau_j.
\end{equation}
In the last part of this section we describe admissible
connections which are compatible with such chirality property. In
contrast to admissibility it is essential to distinguish between
$D$ and~$D^C$, as we will see.

We consider a manifold  of even dimension $2n$ with twisted spinor
bundle $S\oplus S$ and  charge conjugation\footnote{We restrict
ourself to the case  $C\otimes\tau_0=C\oplus C$. The other
possibilities are treated in the same way.} $C\oplus C$. We
suppose that the connection $D$  has an admissible contribution of
the form
\begin{equation}\label{X}
\tfrac{1}{\ell!}F^\ell_{\kappa_1\ldots\kappa_\ell}\gamma^{\kappa_1\ldots\kappa_\ell}\gamma_\mu\otimes\tau_i
+ \tfrac{1}{(2n-\ell)!}
F^{2n-\ell}_{\kappa_1\ldots\kappa_{2n-\ell}}
\gamma^{\kappa_1\ldots\kappa_{2n-\ell}}\gamma^\mu\otimes\tau_j,
\end{equation}
where the two forms are connected by $*F^\ell=w_\ell F^{2n-\ell}$.
We insert this as well as (\ref{duality}) into the connection and
get
\begin{equation*}
\tfrac{1}{\ell!}F^\ell_{\kappa_1\ldots\kappa_\ell}\gamma^{\kappa_1\ldots\kappa_\ell}\gamma_\mu
\big(\mathbbm{1}\otimes \tau_i -w_\ell
(-)^{n}(-)^{\frac{\ell(\ell-1)}{2}}\gamma^*\otimes\tau_j\big).
\end{equation*}

We def\/ine $ \Pi_{ij,w}:=\frac{1}{2}(\mathbbm{1}\otimes \tau_i
+w\gamma^*\otimes\tau_j)$ which has the following properties:
\begin{lemma} \ \ {}

\begin{enumerate}\itemsep=0pt
\item[{\rm 1.}] $\Pi_{ij,w}$ has eigenvalue zero if $(i,j)$ is
none of the pairs $(0,2)$, $(2,0)$, $(1,3)$, or $(3,1)$ and in the
latter cases we have $\Pi_{02,w}^2=\frac{1}{2}w\gamma^*\otimes
\tau_2$ and $\Pi_{13,w}^2=\frac{1}{2}\mathbbm{1}$. \item[{\rm 2.}]
The dimension of the zero eigenspace is $\dim (\ker
\Pi_{ij,w})=\dim S$. \item[{\rm 3.}] For the operators with
eigenvalue zero we have $\Pi_{ij,w}\Pi_{ij,-w}=0$ if $i=j$ or
$ij=0$ and $\Pi_{12,w}\Pi_{12,-w}=\Pi_{03,w}$,
$\Pi_{23,w}\Pi_{23,-w}=\Pi_{01,w}$ but in all cases
\[
{\ker}\,\Pi_{ij,\pm }={\im}\,\Pi_{ij,\mp} .
\]
\end{enumerate}
\end{lemma}
\begin{proof}
The proof is done by taking a look at
\begin{gather*}
\Pi_{ij,w}^2  =\frac{\epsilon_i+\epsilon_j}{4}\mathbbm{1}+\frac{1+\epsilon_{ij}}{4}w\gamma^*\otimes\tau_i\tau_j, \\
\Pi_{ij,w}\Pi_{ij,-w}  =
\frac{\epsilon_i-\epsilon_j}{4}\mathbbm{1}
+\frac{\epsilon_{ij}-1}{4}w\gamma^*\otimes\tau_i\tau_j
\end{gather*}
for the dif\/ferent cases. The kernel of $\Pi_{ij,w}$ which match
with the image of $\Pi_{ij,-w}$ is listed in Table~\ref{kernels}.
\end{proof}

\begin{table}[htb]\caption{The kernels of $\Pi_{ij,w}$ as subsets of $S\oplus S$.}\label{kernels}
$$
{\renewcommand{\arraystretch}{1.5}
\begin{array}{c|c|c|c}
\Pi_{00,w}&\Pi_{01,w}&\Pi_{03,w}&\Pi_{11,w} \\\hline
 S^{-w}\oplus S^{-w}&\big\{ (\eta, -w\gamma^*\eta) | \eta\in S\big\}&S^{-w}\oplus S^{w}& S^{-w}\oplus S^{-w} \\\hline\hline
\Pi_{12,w}&\Pi_{22,w}&\Pi_{23,w}&\Pi_{33,w}\\\hline S^{-w}\oplus
S^{w}&S^{-w}\oplus S^{-w}&\big\{ (\eta, w\gamma^*\eta) | \eta\in
S\big\}&S^{-w}\oplus S^{-w}
\end{array}
}
$$
\end{table}

Due to this lemma we may take $\Pi_{ij,w}$ as a kind of projection
which def\/ines the chirality properties of the spinors $\eta\in
S\oplus S$. (\ref{X}) with $*F^\ell=w_\ell F^{2n-\ell}$
contributes non trivially to an admissible connection in case of a
chirality property of the form $\Pi_{kl,w}\eta=0$ if and only if
${\ker}\,\Pi_{kl,w}\cap{\ker}\,\Pi_{ij,\alpha}\neq 0$ for
$\alpha=w_\ell(-)^n(-)^{\frac{\ell(\ell-1)}{2}}$.

This is the case in (\ref{bsp}) where all projections have the
same image ${\im}\Pi_{11,-}={\im}\Pi_{22,-}= S^+\oplus S^+\subset
S\oplus S$.

We now ask in what way this chirality operator is transferred to
the torsion. The observation  which is summarized in the next
proposition  will, in particular, be used  in
Section~\ref{sectiontorsionfree}.

\begin{proposition}\label{nagative}
Consider a connection $D$ which has a contribution proportional to
a projection cf.~\eqref{projection}. Then the associated part of
the connection $D^C$ as well as the associated part of the torsion
of $D$
 are proportional to the opposite projection.
\end{proposition}

\begin{proof}
We restrict to the the case $i=j=0$ where the connection has a
contribution of the form
$\mathcal{A}_\mu=F_{(\ell)}\gamma^{(\ell)}\gamma_\mu \Pi^\pm $
with $\Pi^\pm=\mathbbm{1}\pm\gamma^*$. The associated part of the
connection $D^C$ is given by
\begin{gather*}
\mathcal{A}^C_\mu
    = F_{(\ell)}( \gamma^{(\ell)}\gamma_\mu \Pi^\pm \big)^C
       = F_{(\ell)}\big(\Delta_1\Delta_\ell\gamma_\mu\gamma^{(\ell)} \mp ( \gamma^{(\ell)}\gamma^* \gamma_\mu )^C\big)\\
\phantom{\mathcal{A}^C_\mu}{}
 = F_{(\ell)}\big( \Delta_1\Delta_\ell\gamma_\mu\gamma^{(\ell)}
 \mp \Delta_1\Delta_{2n-\ell} \gamma_\mu \gamma^{(\ell)}\gamma^*\big)
       = -F_{(\ell)}\gamma_\mu\gamma^{(\ell)}\Pi^\mp,
\end{gather*}
where the last equality is due to the admissibility of the
connection. Furthermore we have
\begin{gather*}
\mathcal{T}_{\mu\nu}
    =\mathcal{A}_\mu\gamma_\nu+\gamma^\nu\mathcal{A}^C_\mu
       =F_{(\ell)}\big(\gamma^{(\ell)}\gamma_\mu \Pi^\pm \gamma_\nu-\gamma_\nu \gamma_\mu\gamma^{(\ell)}\Pi^\mp\big)\\
\phantom{\mathcal{T}_{\mu\nu} }{}
 =F_{(\ell)}\big(\gamma^{(\ell)}\gamma_\mu\gamma_\nu-\gamma_\nu \gamma_\mu\gamma^{(\ell)}\big)\Pi^\mp
       =F_{(\ell)}\big(\gamma^{(\ell)}\gamma_{\mu\nu} +\gamma_{\mu\nu}\gamma^{(\ell)}\big)\Pi^\mp .
\end{gather*}
The proof for
$\mathcal{A}_\mu=F_{(\ell)}\gamma_\mu\gamma^{(\ell)}\Pi^{\pm}$ or
$(i,j)\neq(0,0)$ is almost the same.
\end{proof}

\subsection{Jacobi versus Bianchi}\label{JacobiBianchi}

In this section we consider a graded manifold of the form $\hat
M=(M,\Lambda\Gamma S)$ and calculate commutators of the vector
f\/ields $\imath(\phi)$ which have been def\/ined in the
preliminaries. The (graded) Jacobi identity on the (super) Lie
algebra of vector f\/ields will be seen to be related to the
Bianchi identities.

We recall the inclusions of the vector f\/ields on $M$, the
spinors, and of the endomorphisms of~$S$  into the vector f\/ields
on $\hat M$ as given in Section~\ref{inclusions}. Due to the fact
that we will f\/ix a connection~$D$ on~$S$, we will drop the index
and will write $\jmath:\Gamma S\oplus\X(M)\to\X(\hat M)$ for the
inclusions.
\begin{proposition}\label{propoo}
We consider  the  graded manifold $(M,\Lambda \Gamma S)$ and a
connection $D$ on $S$ which defines the inclusion $\jmath$ and the
map $\imath:\Gamma S\to \X(\hat M)$. Furthermore we consider  a
linear subspace
 $\mathcal{K}\subset\big\{\phi\in\Gamma S\,|\,D^C\eta=0\big\}$
 such that $(\mathcal{K},D)$ is admissible. Then  the following holds for all $\varphi,\psi\in \mathcal{K}$
\begin{gather}
\big[\imath(\varphi),\imath(\psi)\big]
    =  \mathfrak{B}(R;\varphi,\psi) +\tfrac{1}{2}\mathfrak{D}(\mathcal{T};\varphi,\psi),
\label{oo}
\end{gather}
where we use the short notations
\begin{gather}
\mathfrak{B}(R;\varphi,\psi) =\gamma^\mu\varphi\wedge\gamma^\nu\psi\owedge R_{\mu\nu},  \label{B}\\
\mathfrak{D}(\mathcal{T};\varphi,\psi) =
       \big( \gamma^\mu\varphi\wedge  \mathcal{T}_{\mu\nu}\psi
    +\gamma^\mu\psi\wedge \mathcal{T}_{\mu\nu} \varphi\big) \otimes D^\nu.  \label{D}.
\end{gather}
\end{proposition}

\begin{proof}
For all $\phi,\psi\in\Gamma S$ we have
\begin{gather*}
 \big[\imath(\phi),\imath(\psi)\big]
    = \gamma^\mu\phi\owedge \big[\jmath(e_\mu),\gamma^\nu\psi\otimes \jmath(e_\nu)\big]
        +\gamma^\nu\psi\wedge D_\nu(\gamma_\mu\phi)\otimes \jmath(e^\nu)\\
\phantom{\big[\imath(\phi),\imath(\psi)\big]}{}    =
\gamma^\mu\phi\wedge D_\mu(\gamma^\nu\psi)\otimes \jmath(e_\nu)
        + \gamma^\mu\phi\wedge \gamma^\nu\psi\owedge \big[\jmath(e_\nu),\jmath(e_\nu)\big]\\
\phantom{\big[\imath(\phi),\imath(\psi)\big]}{}     =\
\mathfrak{B}(R;\phi,\psi)
        +\gamma^\mu\phi\wedge  \hat D_\mu \gamma_\nu\,\psi\otimes D^\nu
        +\gamma^\mu\psi\wedge \hat D_\mu \gamma_\nu\, \phi\otimes D^\nu\\
\phantom{\big[\imath(\phi),\imath(\psi)\big]=}{}
+\gamma^\mu\phi\wedge \gamma^\nu D^C_\mu\psi \otimes D_\nu
            +\gamma^\mu\psi\wedge \gamma^\nu D^C_\mu\phi \otimes D_\nu.
\end{gather*}
In particular, these relations reduce to (\ref{oo}) if we restrict
to $\mathcal{K}$.
\end{proof}
\begin{corollary}
Consider an admissible metric connection on $S$, i.e.\ with skew
symmetric torsion $T_{\mu\nu\kappa}=2A_{\mu\nu\kappa}$. In this
case \eqref{oo}  is given by
\begin{gather*}
\big[\imath(\varphi),\jmath(\psi)\big]   =
    \tfrac{1}{4}R_{\mu\nu\kappa\rho}\gamma^\mu\varphi\wedge\gamma^\nu\psi\owedge\gamma^{\kappa\rho}
    + \tfrac{1}{2}T_{\mu\nu\kappa}\gamma^\mu\varphi\wedge\gamma^\nu\psi\otimes D^\kappa.
\end{gather*}
\end{corollary}

For the following calculations we restrict to the case that the
spinors belong to an admissible subspace
$\mathcal{K}\subseteq\big\{\eta\in\Gamma S\,|\,D^C\eta=0\big\}$
\begin{gather}
 \big[\imath(\varphi),\mathfrak{B}(R;\eta,\xi) \big]
    = \tfrac{1}{2}\gamma_\kappa\varphi\wedge \mathcal{T}^{\kappa\mu}\eta\wedge\gamma^\nu\xi\owedge R_{\mu\nu}
                    + \tfrac{1}{2}\gamma_\kappa\varphi\wedge \gamma^\mu\eta \wedge
                    \mathcal{T}^{\kappa\nu}\xi\owedge R_{\mu\nu}        \nonumber\\
\phantom{\big[\imath(\varphi),\mathfrak{B}(R;\eta,\xi) \big]= }{}
-\gamma^\mu\eta\wedge\gamma^\nu\xi\wedge
ad^C_{R_{\mu\nu}}\gamma^\kappa\varphi\otimes D_\kappa
        +  \gamma^\kappa\varphi\wedge\gamma^\mu\eta\wedge\gamma^\nu\xi\owedge (D_\kappa R)_{\mu\nu},\label{d}\\
\big[\imath(\varphi) , \mathfrak{D}(\mathcal{T};\eta,\xi)\big]
     = \tfrac{1}{2}\gamma^\kappa\varphi\wedge\mathcal{T}_{\kappa\mu}\eta\wedge\mathcal{T}^{\mu\nu}\xi\otimes D_\nu
            + \tfrac{1}{2}\gamma^\kappa\varphi\wedge\mathcal{T}_{\kappa\mu}\xi\wedge\mathcal{T}^{\mu\nu}\eta\otimes D_\nu
        \nonumber \\
\phantom{\big[\imath(\varphi) ,
\mathfrak{D}(\mathcal{T};\eta,\xi)\big] =}{} +
\tfrac{1}{2}\mathcal{T}^{\nu\kappa}\varphi\wedge\gamma^\mu\eta\wedge\mathcal{T}_{\mu\kappa}\xi\otimes
D_\nu
                    + \tfrac{1}{2}\mathcal{T}^{\nu\kappa}\varphi\wedge\gamma^\mu\xi\wedge\mathcal{T}_{\mu\kappa}\eta\otimes D_\nu
        \nonumber \\
\phantom{\big[\imath(\varphi) ,
\mathfrak{D}(\mathcal{T};\eta,\xi)\big] =}{}  +
\gamma^\kappa\varphi\wedge\gamma^\mu\eta\wedge \hat
D_\kappa\mathcal{T}_{\mu\nu}\xi\otimes D^\nu
                    + \gamma^\kappa\varphi\wedge\gamma^\mu\xi\wedge \hat D_\kappa\mathcal{T}_{\mu\nu}\eta\otimes D^\nu \nonumber\\
\phantom{\big[\imath(\varphi) ,
\mathfrak{D}(\mathcal{T};\eta,\xi)\big] =}{}      +
\gamma^\kappa\varphi\wedge\gamma_\mu\eta\wedge\mathcal{T}^{\mu\nu}\xi\owedge
R_{\kappa\nu}
                    + \gamma^\kappa\varphi\wedge\gamma_\mu\xi\wedge\mathcal{T}^{\mu\nu}\eta\owedge R_{\kappa\nu}.\label{e}
\end{gather}

From (\ref{d}) and (\ref{e}) we read of the terms of dif\/ferent
order in $
\big[\imath(\varphi),\big[\imath(\eta),\imath(\xi)\big]\big]$:
\begin{gather}
 \big[\imath(\varphi),\big[\imath(\eta),\imath(\xi)\big]\big]^{(3,0)}
     =\tfrac{1}{4}\big( \gamma^\kappa\varphi\wedge\mathcal{T}_{\kappa\mu}\eta\wedge\mathcal{T}^{\mu\nu}\xi
                +\gamma^\kappa\varphi\wedge\mathcal{T}_{\kappa\mu}\xi\wedge\mathcal{T}^{\mu\nu}\eta \nonumber\\
\qquad\qquad  {}+
\mathcal{T}^{\nu\kappa}\varphi\wedge\gamma^\mu\eta\wedge\mathcal{T}_{\mu\kappa}\xi
                +\mathcal{T}^{\nu\kappa}\varphi\wedge\gamma^\mu\xi\wedge\mathcal{T}_{\mu\kappa}\eta  \big)
        \otimes D_\nu \nonumber\\
\qquad   \qquad  {} +\Big(
\tfrac{1}{2}\gamma^\kappa\varphi\wedge\gamma^\mu\eta\wedge \hat
D_\kappa\mathcal{T}_{\mu\nu}\xi
                + \tfrac{1}{2}\gamma^\kappa\varphi\wedge\gamma^\mu\xi\wedge \hat D_\kappa\mathcal{T}_{\mu\nu}\eta \nonumber\\
\qquad   \qquad  {}- \gamma^\kappa\eta\wedge\gamma^\mu\xi\wedge
ad^C_{R_{\kappa\mu}}\gamma_\nu\varphi
            \Big) \otimes D^\nu,\label{30} \\
 \big[\imath(\varphi),\big[\imath(\eta),\imath(\xi)\big]\big]^{(4,1)}
    =\tfrac{1}{2}\big(  \gamma^\mu\varphi\wedge\gamma_\kappa\eta\wedge\mathcal{T}^{\kappa\nu}\xi
                    + \gamma^\mu\varphi\wedge\gamma_\kappa\xi\wedge\mathcal{T}^{\kappa\nu}\eta \nonumber \\
\qquad   \qquad  {}+\gamma_\kappa\varphi\wedge\gamma^\mu\xi \wedge
\mathcal{T}^{\kappa\nu}\eta
                + \gamma_\kappa\varphi\wedge \gamma^\mu\eta \wedge \mathcal{T}^{\kappa\nu}\xi \big)
        \owedge R_{\mu\nu} \nonumber\\
\qquad   \qquad  {}
+\gamma^\kappa\varphi\wedge\gamma^\mu\eta\wedge\gamma^\nu\xi\owedge
(D_\kappa R)_{\mu\nu}.\label{41}
\end{gather}

The Jacobi identity, i.e.\ the vanishing of
 $ \cycl\limits_{\varphi,\eta,\xi}\big[\imath(\varphi),\big[\imath(\eta),\imath(\xi)\big]\big]$
 holds independently for  the terms of dif\/ferent degree -- here $\cycl$ denotes the graded cyclic sum.
More precisely:
\begin{itemize}\itemsep=0pt
\item The cyclic sums of the f\/irst summands  in (\ref{30}) and
(\ref{41}) vanish due to the symmetry of the involved objects.
\item The vanishing of the  cyclic sum of the second summand in
(\ref{41}) is equivalent to  the Bianchi identity
(\ref{BianchiDR}). \item The cyclic sum of the second summand in
(\ref{30}) vanishes due to the algebraic Bianchi-identity of the
curvature of the Levi-Civita connection. This is due to the
following supplement to Proposition~\ref{Bianchi}.
\end{itemize}

\begin{lemma}\label{Bianchi2}
Let $\mathcal{K}\subset\big\{\eta\in\Gamma S\,|\,D^C\eta=0 \big\}$
such that  $(D,\mathcal{K})$ is admissible. Let $\mathcal{T}$ be
the torsion of  $D$. Then \eqref{BianchiDT} in Proposition~{\rm
\ref{Bianchi}} reduces to
\begin{equation*}
\big(\hat D_{[\kappa}\mathcal{T}_{\mu]\nu}-{\rm
ad}^C_{R_{\kappa\mu}} \gamma_\nu-
R^0_{\kappa\mu\nu\lambda}\gamma^\lambda\big)\eta=0
\end{equation*}
for all $\eta\in\mathcal{K}$.
\end{lemma}
This yields
\begin{corollary}
Let $\mathcal{K}$, $D$ and $\mathcal{T}$ as before. For all
spinors $\varphi,\eta,\xi\in\mathcal{K}$ the following holds
\begin{equation*}
\cycl_{\varphi,\eta,\xi} \Big\{ \big(\hat
D_{[\kappa}\mathcal{T}_{\mu]\nu}-{\rm
ad}^C_{R_{\kappa\mu}}\gamma_\nu\big)\xi \wedge
\gamma^\kappa\varphi\wedge\gamma^\mu\eta \Big\} = 0 .
\end{equation*}
\end{corollary}

\begin{remark}
As we saw above, the action of $\hat
D_{[\kappa}\mathcal{T}_{\mu]\nu}-{\rm
ad}^C_{R_{\kappa\mu}}\gamma_\nu$ on $\mathcal{K}$ coincides with
the action of the curvature of the Levi-Civita connection $R^0$ on
$\mathcal{K}$. If $D$ is admissible this yields a way to express
$R^0$ in terms of  $R$ and $\mathcal{T}$. Let furthermore  $D$ be
metric, i.e.\ a connection with totally skew symmetric torsion.
Then the above expression may be written as
\[
R^0_{\kappa\lambda\mu\nu}=R_{\kappa\lambda\mu\nu}-
D_{[\kappa}^TT_{\lambda]\mu\nu}-\tfrac{1}{4}T_{\kappa\lambda\rho}T_{\mu\nu}{}^\rho
- \sigma^T_{\kappa\lambda\mu\nu}
\]
with $\sigma^T_{\kappa\lambda\mu\nu}= 3 T_{\rho[\kappa\lambda
}T_{\mu]\nu}{}^\rho$ which is indeed a 4-form. This is due to
\cite{FriedIvanov3} or \cite{IvanovPapado}. Here $D^T$ denotes the
connection which dif\/fers from $D$ by
\[
(D_X-D^T_X)T(Y,Z)
=\tfrac{1}{2}T(T(X,Y),Z)+\tfrac{1}{2}T(Y,T(X,Z)),
\]
i.e.\
$(D_\mu-D^T_\mu)T_{\kappa\lambda\nu}=T_{\rho\nu[\lambda}T_{\kappa]\mu}{}^{\rho}$.
\end{remark}

\section{Applications and examples}\label{exmps}

\subsection{Torsion freeness}\label{sectiontorsionfree}

We consider a connection $D$ on the spinor bundle $S$ and
$\mathcal{K}\subseteq\Gamma S$ such that $(D,\mathcal{K})$ is
admissible. In (\ref{oo}) we def\/ined the map
$\mathfrak{D}:S^2(\Gamma S)\to \Lambda^2\Gamma S\otimes \X(M)$
which motivates the following def\/inition.
\begin{definition}\label{definitiontorsionfree}
Let $D$ be a connection on $S$ with torsion $\mathcal{T}$ and
$\mathcal{K}\subseteq\Gamma S$ such that $(D,\mathcal{K})$ is
admissible.
\begin{enumerate}\itemsep=0pt
\item We call $(D,\mathcal{K})$  {\sl torsion free}  if
$\mathfrak{D}(\mathcal{T};\eta,\xi)=0\  \text{ for all }\
\eta,\xi\in\mathcal{K}$. \item We call $(D,\mathcal{K})$  {\sl
strongly torsion free} if $\mathcal{T}_{\mu\nu}\eta=0\ \text{ for
all }\ \eta\in\mathcal{K}$.
\end{enumerate}
And in view of (\ref{oo})
\begin{enumerate}\itemsep=0pt
\addtocounter{enumi}{2}\item We call  $(D,\mathcal{K})$  {\sl
flat} if $\mathfrak{B}(R;\varphi,\psi)
=\mathfrak{D}(\mathcal{T};\varphi,\psi)=0$ for all
$\varphi,\psi\in\mathcal{K}$.
\end{enumerate}
\end{definition}

There are two natural problems:  f\/irstly f\/ix $D$ and restrict
 $\mathcal{K}$ such that one of the properties are obtained,
 secondly look for conditions on the connection -- or the torsion --
 such that an admissible set $\mathcal{K}$ is ``as large as needed''.

Of course,  admissible subsets
$\mathcal{K}\subseteq\big\{\eta\in\Gamma S\,|\,D^C\eta=0\big\}$
will be of particular interest. Due to Theorem \ref{KillingA} the
Killing vector f\/ields which we obtain by
$\{\mathcal{K},\mathcal{K}\}$ are parallel with respect to the
Levi-Civita connection if $\mathcal{K}$  is strongly torsion free.
Therefore, to get  non parallel Killing vector f\/ields by $C_1$,
it is necessary for the connection $D$ on $S$ to admit a part
which contribution to the torsion acts non trivially on
$\mathcal{K}$.

\subsubsection[On strongly torsion freeness in $R^n$]{On strongly torsion freeness in $\boldsymbol{\r^n}$}

We consider f\/lat $\r^n$ with spinor bundle $S$ and connection
\begin{equation*}
D^C_X\psi=X(\psi)-\mathcal{A}^C_\mu\psi.
\end{equation*}
where the potential $\vec
A=(\mathcal{A}^C_1,\ldots,\mathcal{A}^C_n)$ is constructed from
forms on $\r^n$ with constant coef\/f\/i\-cients.

\begin{example}\label{exsun}
Consider $\r^{2n}$ with connection
$D^C_\mu=d_\mu-\mathcal{A}^C_\mu$ on its spinor bundle. Let
$\mathcal{A}^C$  be determined by a three-form $F$, moreover $F$
shall be a one-form with values in $\mathfrak{su}(n)$. Then
$\mathfrak{hol}\subset\mathfrak{su}(n)$ and there exist two
parallel pure spinors $\eta,\bar\eta$  which are associated via
charge conjugation. These spinors obey
$\mathfrak{B}(R;\eta,\eta)=\mathfrak{B}(R,\bar\eta,\bar\eta)=0$.
We use the decomposition
$\c^{2n}=\mathbf{n}\oplus\bar{\mathbf{n}}$ where the complex
structure obeys  $\mathbf{n}\bar\eta=\bar{\mathbf{n}}\eta=0$. If
$F\in\Lambda^3\c^{2n}\cap(\bar{\mathbf{n}}\otimes\mathfrak{su}(n))$
the torsion acts trivially on $\eta$. In this case the subspace
spanned by this sole odd generator is strongly torsion free, in
particular, it would have vanishing center. If
$F\in\Lambda^3\c^{2n}\cap(\mathbf{n}\otimes\mathfrak{su}(n))$ the
same holds for $\bar\eta$. We emphasize that in both cases the
three-form is not real and that for a real three-form a trivial
action on one of the spinors is only possible in case of vanishing
torsion.
\end{example}

\begin{example}
Suppose $\mathcal{A}$ is obtained by a constant form and
$\mathcal{A}_X\propto \Pi^+$ ($\Pi^-$) for a projection~$\Pi^\pm$
cf.~(\ref{projection}). Due to Proposition~\ref{nagative} $D^C$
and $\mathcal{T}$ are proportional to the opposite
projection~$\Pi^-$~($\Pi^+$). So  $\mathcal{K}$ spanned by the
constant positive (negative) spinors is strongly torsion free.
\end{example}

The last example can be generalized to
\begin{remark}
Strongly torsion freeness can not be achieved by pure chirality
considerations due to Proposition \ref{nagative}, when we want to
deal with spinors which are not Levi-Civita parallel. In this case
strongly torsion freeness leads to new algebraic constraints on
the f\/ields.
\end{remark}

We will discuss torsion free structures which are not strongly
torsion free in Section~\ref{examplePure} (generalizing Example
\ref{exsun})  and \ref{examplebranes}.

\subsubsection[A comment on differentials]{A comment on dif\/ferentials}
As we mentioned in the introduction and as performed in
\cite{Papado2} we may take the vector f\/ield
$\imath_D(\eta)=\gamma^\mu\eta\otimes \jmath_D(e_\mu)$ as
degree-one operator on $\Lambda S$ and look for conditions such
that this operator is a dif\/ferential. We immediately get
\begin{proposition}\label{differential}
Let $D$ be a connection on a spinor bundle $S$ over the (pseu\-do)
Riemann\-ian manifold $M$. Consider the vector field
$\imath(\eta)$ on the graded manifold $(M,\Gamma\Lambda S)$. Let
$D^C\eta=0$, then $\imath(\eta)$ is a dif\/ferential on
$\Gamma\Lambda S$ if and only if $(D,\{\eta\})$ is flat.
\end{proposition}

When we consider admissible subspaces  $\mathcal{K}$ of order one
we have to  take the collection of all elements in
$\big\{\mathfrak{B}(R;\eta,\xi) |\,D^C\eta=D^C\xi=0\big\}$ and
$\big\{\mathfrak{D}(\mathcal{T};\eta,\xi)
|\,D^C\eta=D^C\xi=0\big\}$ and discuss whether or not these terms
vanish. In particular, if the dimension of $\mathcal{K}$ is large
the conditions on the torsion are very restrictive. When we
consider the  dif\/ferential point of view we only have to discuss
the terms $\mathfrak{B}(R;\eta,\eta)$ and
$\mathfrak{D}(\mathcal{T};\eta,\eta)$ for one f\/ixed spinorial
entry.

In \cite{Papado2} and \cite{Klinker4} the condition on
$\mathfrak{B}$  is discussed for the untwisted case. The twisted
case is touched when the authors discuss real spinors. The main
emphasis  is on metric connections $D$ of  holonomy $\mathfrak{g}
\subset\mathfrak{so}(n)\subset\mathfrak{sl}(2^{[\frac{n}{2}]})$
with $\mathfrak{g}=\mathfrak{su}(\frac{n}{2})$,
$\mathfrak{sp}(\frac{n}{4})$, $\mathfrak{spin}(7)$ if $n=8$, or
$\mathfrak{g}_2$ if $n=7$. The discussion in \cite{Klinker4} is
restricted to the torsion free Levi-Civita connection. If we want
to cover non-torsion free metric connections -- or general
spinorial connections -- we have to take into account the
$\mathfrak{D}$-contribution which yields further restrictions and
we recall Example \ref{exsun} and the examples below.


\subsection{Parallel pure spinors}\label{examplePure}

We consider a Riemannian manifold $M$ of even dimension $2n\geq
4$. Consider a pure spinor $\eta\in\Gamma S$. We will discuss
conditions on a connection $D$ such that
$\mathfrak{B}(R;\eta,\eta)$ or
$\mathfrak{D}(\mathcal{T};\eta,\eta)$ vanish. As before, the case
of a $D^C$-parallel pure spinor is of particular interest due to
Theorem~\ref{KillingA}, Section~\ref{JacobiBianchi}, and
Proposition~\ref{differential}. Although we deal with forms of
arbitrary degree, we always specialize  to the metric case.

A pure spinor is characterized by the following two equivalent
conditions (compare \cite{Chevalley,KathPure}).
\begin{enumerate}\itemsep=0pt
\item[(1)] The space $\{X\in TM\ |\ X\eta=0\}$ has maximal
dimension, namely $n$.\label{pure1} \item[(2)] $C_k(\eta,\eta)=0$
for all $k\neq n$.\label{pure2}
\end{enumerate}
Furthermore a pure spinor is either of positive or of negative
chirality and  the vector f\/ield $\{\eta,\eta\}$ vanishes. The
symmetry $\Delta_k$ and the chirality of $C_k$ are given by
$$
 \begin{array}{c||c|c|c|c}
2n\,{\text{mod}}\,8 &0      &2      &4      &6      \\\hline\hline
 \Delta_{2m}    &(-)^m      &\pm(-)^m   &-(-)^m     &\mp(-)^m   \\\hline
 \Delta_{2m+1}  &\pm(-)^m   &(-)^m      &\mp(-)^m   &-(-)^m     \\\hline\hline
 {\text{chirality}} &{\text{non chiral}}    &\text{chiral}  &\text{non chiral}  &\text{chiral}
\end{array}
$$

The dif\/ferent signs belong to the choice of charge conjugation.
Chiral means $C: S^\pm\otimes S^\mp\to \c$ and non-chiral (nc)
means $C:S^\pm\otimes S^\pm\to \c$. Examining this table yields
that the second part~2 in the characterization may be relaxed as
follows
\begin{enumerate}\itemsep=0pt
\item[$(2')$]The chiral (or anti-chiral) spinor $\eta$ is pure if
$C_k(\eta,\eta)=0$ for all $k-n\equiv 0\,{\text{mod}}\,4$, $k\neq
n$.\label{pure3}
\end{enumerate}
In particular $C_n$ has symmetry $\Delta_n=1$ in all cases.

We  take a closer look at $\mathfrak{B}(R;\eta,\eta)
=\gamma^\mu\eta\wedge\gamma^\nu\eta\owedge R_{\mu\nu}$. We use the
Fierz identity (\ref{Fierz}) to to rewrite this expression.
\begin{gather*}
\gamma^{[\mu}\eta\wedge\gamma^{\nu]}\eta  =
   \frac{1}{\dim S}\sum_{\Delta_k=-1}\frac{\Delta_0(\Delta_0\Delta_1)^k}{ k!}
    C(\gamma^{[\mu}\varphi,\gamma_{(k)}\gamma^{\nu]}\psi) \gamma^{(k)} \\
\phantom{\gamma^{[\mu}\eta\wedge\gamma^{\nu]}\eta}{} =
\frac{1}{\dim S}\sum_{\Delta_k =-1}
    (-\Delta_0\Delta_1)^{k+1} \Big( \frac{1}{k!}C(\gamma^{\mu\nu (k) }\psi,\varphi) \gamma_{(k)}\\
\phantom{\gamma^{[\mu}\eta\wedge\gamma^{\nu]}\eta=}{}
  + \frac{1}{(k-2)!}C(\gamma_{(k-2) }\psi,\varphi) \gamma^{\mu\nu (k-2)}\Big)\\
\phantom{\gamma^{[\mu}\eta\wedge\gamma^{\nu]}\eta}{}=
 \frac{(-\Delta_0\Delta_1)^{n+1}}{\dim S}\Big( \frac{1}{(n-2)!}C(\gamma^{\mu\nu(n-2)}\eta,\eta)\gamma_{(n-2)}
 +\frac{1}{n!}C(\gamma_{(n)}\eta,\eta)\gamma^{\mu\nu(n)} \Big).
\end{gather*}
The second last equality holds because of (\ref{iden}) and the
last due to the fact that only the summands with $k=n-2$ and
$k=n+2$ survive. Furthermore we needed
$1=\Delta_n=-\Delta_{n-2}=-\Delta_{n+2}$. Using the duality
relation (\ref{duality}) to manipulate the f\/irst or second
summand, we get the following two equivalent expressions
\begin{gather}
\gamma^{[\mu}\eta\wedge\gamma^{\nu]}\eta =
\frac{(-\Delta_0\Delta_1)^{n+1}}{n! \dim S}
 C(\gamma_{(n)}\eta,\eta) \gamma^{\mu\nu(n)} (\mathbbm{1}-(-)^nw_\eta\gamma^*)\nonumber \\
\intertext{and} \gamma^{[\mu}\eta\wedge\gamma^{\nu]}\eta =
\frac{(\Delta_0\Delta_1)^{n+1}}{(n-2)! \dim S}
 C(\gamma^{\mu\nu(n-2)}\eta,\eta) \gamma_{(n-2)} (\mathbbm{1}-(-)^n w_\eta\gamma^*), \label{selfdual}
\end{gather}
where $w_\eta$ is def\/ined by $\gamma^*\eta=w_\eta\eta$.

Suppose $\dim M= 4$. Then (\ref{selfdual}) is self dual if $\eta$
is of negative chirality and anti-self dual if $\eta$ is of
positive chirality in the sense that
\[
\tfrac{1}{2}\epsilon_{\rho\sigma\mu\nu}\gamma^{[\mu}\eta\wedge\gamma^{\nu]}\eta
=
  -w_\eta \gamma_{[\rho}\eta\wedge\gamma_{\sigma]}\eta.
\]
This yields
\begin{proposition} Let $M$ be of dimension four and the pure spinor $\eta$ be of negative (positive) chirality.
Then $\mathfrak{B}(R;\eta,\eta)$ vanishes if the curvature $R$ of
$D$ is self dual (resp.\ anti-self dual).
\end{proposition}

The last proposition  is an extension of the result we obtained
in \cite{Klinker4} where we examined the four dimensional case
with $D=\nabla$ and holonomy $\mathfrak{su}(2)$ which implies
self-duality of the curvature tensor $R^0$. Moreover in dimension
four there is a further symmetry which yields
$\mathfrak{B}(R^0;\eta,\eta^C) =0$ for  the parallel pure spinors
$\eta$ and its parallel pure charge conjugated $\eta^C$.

Self duality of the curvature tensor as a necessary condition for
the vanishing of $\mathfrak{B}(R;\eta,\eta)$  is too restrictive.
Suppose $\eta$ is positive so that
$\gamma_{[\mu}\eta\wedge\gamma_{\nu]}\eta$ is anti-self dual. This
is half of the game. More precisely we f\/ind
$\gamma_{[\mu}\eta\wedge\gamma_{\nu]}\eta$ in the $\Lambda^{2,0}$
part of anti-self dual forms $\Lambda^2_-\otimes\c$. Here
$\Lambda^{2,0}$ is def\/ined by  the complex structure given by
$\eta$ (compare~\cite{LawMich}). If we use complex matrices
$\{\gamma^a,\gamma^{\bar a}\}_{1\leq a,\bar a \leq 2}$
 associated to this complex structure, i.e.\ $\gamma^{\bar a}\eta=0$,
 and write $R$ in this frame as $R_{ab},R_{a\bar b},R_{\bar a\bar b}$
 the necessary condition for the vanishing of $\mathfrak{B}(R;\eta,\eta)$ is $R_{12}=0$.

If the connection $D$, and so the curvature $R$,  is in a real
representation the vanishing of the $\Lambda^{2,0}$-part of the
curvature is equivalent to two of the three self duality
equations. Furthermore we have\footnote{$A$ denotes the matrix
which def\/ines the charge conjugation $\varphi^C:=A\varphi^*$
compare \cite{Klinker4}.} $R_{ij}=A R^*_{\bar\imath\bar\jmath}
A^{-1}$  and the $\Lambda^{0,2}$ part $R_{\bar 1\bar 2}$ vanishes,
too. So the condition for the vanishing of $\mathfrak{B}$ reduces
to  $R\in\Lambda^{1,1}$. The part  which prevent the curvature
from being self dual is the trace of the $\Lambda^{1,1}$-part.
This is due to the isomorphism $\Lambda^2_+=\Lambda^{1,1}_0$,
cf.~\cite{AtHiSi}.

Similar considerations as in the four dimensional case can be made
for arbitrary even dimension. For this we introduce complex
coordinates associated to the null space of $\eta$,
$\{\gamma^a,\gamma^{\bar a}\}_{1\leq a,\bar a \leq n}$ with
$\gamma^{\bar a}\eta=0$ . The only surviving part of the form
which is associated to $\eta$ via the Fierz identity is
$C(\gamma^{1\ldots n}\eta,\eta) \gamma_{1\ldots n}$ with only
unbarred indices. So (\ref{selfdual}) reads as
\begin{equation*}
\gamma^{\mu}\eta\wedge\gamma^{\nu}\eta \owedge R_{\mu\nu} =
        \eta^{(n)} \varepsilon^{a_1\ldots a_n}\gamma_{a_1\ldots a_{n-2}}
        (\mathbbm{1}-(-)^nw_\eta\gamma^*) \owedge R_{a_{n-1}a_n}
\end{equation*}
with $\eta^{(n)}:=\frac{(-\Delta_0\Delta_1)^{n+1}}{(n-2)! \dim S}
C(\gamma^{1 \ldots n}\eta,\eta)$ and $\varepsilon^{a_1\ldots a_n}$
the totally skew-symmetric symbol of unbarred indices. This yields
\begin{proposition}\label{proppurecurv}
Let $\eta$ be a pure spinor on the even dimensional manifold $M$.
Then $\mathfrak{B}(R;\eta,\eta)$ vanishes if and only if
\begin{equation}
\varepsilon_{a_1\ldots a_n}\gamma^{a_1\ldots a_{n-2}}
(\mathbbm{1}-(-)^n w_\eta\gamma^*) \owedge R^{a_{n-1}a_n}=0.
\label{nessuffB}
\end{equation}
Here the sum is over the unbarred indices with respect to the
complex structure given by the pure spinor $\eta$.
\end{proposition}

A class of connections for which the above is applicable is given
in the following corollary. The proof needs the decomposition of
$\Lambda^2$ which can be taken from the discussion of the four
dimensional case.
\begin{corollary}
Let $D$ be a metric connection on $M$, and suppose it is of
holonomy $\mathfrak{su}(n)$. Then condition \eqref{nessuffB} holds
for the two parallel pure spinors.
\end{corollary}

Using the complex coordinates which have been introduced above,
the condition $R_{ab}=0$ as a necessary condition for
$\mathfrak{B}(R;\eta,\eta)=0$ could be seen directly from
(\ref{B}). Nevertheless, we used the Fierz identity here to draw a
connection to the forms def\/ined by the spinor $\eta$ and to make
the condition more precise.

We turn to the torsion dependent term
$\mathfrak{D}(\mathcal{T};\eta,\eta)$ and distinguish the two
cases
$\mathcal{T}_{\mu\nu}=\frac{1}{(\ell-2)!}F_{\mu\nu(\ell-2)}\gamma^{(\ell-2)}$
and
$\mathcal{T}_{\mu\nu}=\frac{1}{\ell!}F^{(\ell)}\gamma_{\mu\nu(\ell)}$.
In both cases we use the Fierz identity as well as (\ref{iden})
and condition ($2'$) above  and get after some careful
calculations
\begin{gather}
 \frac{1}{(\ell-2)!}  \gamma^\mu  \eta\wedge F_{\mu\nu(\ell-2)}\gamma^{(\ell-2)}\eta \nonumber\\
\qquad{}= \frac{1}{\dim
S}\sum_{\Delta_k=-1}(\Delta_0\Delta_1)^{k+1}
        \frac{ (-)^{ \frac{m(m-2k-1)}{2}}(-)^{k-1}(n-k)}{m!(k-m)!(\ell-m-1)!}\times \nonumber \\
\qquad\phantom{=}{}
  \times F_{\nu(m)}{}^{(\ell-1-m)}  C\big( \gamma_{(k-m)(\ell-1-m)} \eta,\eta\big) \gamma^{(m)(k-m)}\label{eins} \\
\intertext{and}
 \frac{1}{\ell!}  \gamma^\mu  \eta\wedge F^{(\ell)} \gamma_{\mu\nu(\ell)}\eta
= \frac{1}{\dim S}\sum_{\Delta_k=-1}(\Delta_0\Delta_1)^{k+1}
       \frac{(-)^{ \frac{ m(m-2k-1)}{2}}(-)^{m+1}(n-k)}{( m+1)!(k-m-1)!(\ell-m)!}\times \nonumber \\
\qquad{} \times \Big(  (m+1) F^{(m)}{}_{(\ell-m)}
    C\big( \gamma^{(k-m-1)(\ell-m)}\eta,\eta\big)\gamma_{\nu(m)(k-m-1)}\nonumber \\
\qquad{} +(-)^\ell (\ell-m) F_{(m+1)}{}^{(\ell-m-1)}
    C\big( \gamma_{(k-m-1)(\ell-m-1) \nu} \eta,\eta\big) \gamma^{(m+1)(k-m-1)}\Big)\label{zwei}
\end{gather}
with $m=\tfrac{1}{2}(k+\ell-n-1)$. This may be used to get
conditions on the forms and their contribution to the connection
$D$ to let $\mathfrak{D}(\mathcal{T};\eta,\eta)$ vanish. We will
not explicitly use this formulas in the next example, but we will
see that this would have been possible.

\begin{example}
We turn again to  the case of  dimension four. In the case
$\ell=3$, i.e.\ the case of metric connection
$D_\mu-\nabla_\mu^0=D_\mu^C-\nabla_\mu^0=\mathcal{A}_\mu=T_{\mu\nu\kappa}\gamma^{\nu\kappa}$
the term
\[
\gamma^\mu\eta\wedge\mathcal{T}_{\mu\nu}\eta =
T_{\mu\nu\kappa}\gamma^\mu\eta\wedge\gamma^\kappa\eta
\]
vanishes in the case of  self duality. We recall the decomposition
$\Lambda^2\otimes\Lambda^1
=\Lambda^1\oplus\Lambda^3\oplus\Lambda^{(2,1)}$. If we denote the
projections on $\Lambda^1\simeq\Lambda^3$ and $\Lambda^3$ by
$\pi_1$ and $\pi_3$  respectively, we have
\begin{equation}
T\in \Lambda^2_\pm\otimes\Lambda^1\ \Longleftrightarrow\ *
\pi_1(T)= \pm \pi_3(T).\label{dualforms}
\end{equation}

This example f\/its into the discussion of admissible connections,
in particular, when we added ``non-allowed'' forms to the
connection in the artif\/icial way  (\ref{starcontributioneven}).
Moreover if we would have taken an arbitrary one-form $V^\kappa
\gamma_{\mu\kappa}\gamma^*$ and three-form
$T_{\mu\nu\kappa}\gamma^{\nu\kappa}$ as contributions to
$\mathcal{A}=D-\nabla$, equations (\ref{eins}) and (\ref{zwei})
would have yield exactly the right hand side of (\ref{dualforms}).
\end{example}

As before we may generalize the result to dimensions greater than
four. When we consider three-form potentials we see  that the
$\mathfrak{D}$- and the $\mathfrak{B}$-term have similar shape. So
we get
\begin{proposition}\label{proppuretor}
Let $\eta$ be a pure spinor  and $D$ be constructed by a $3$-form.
Then $\mathfrak{D}(\mathcal{T};\eta,\eta)=0$ if
\begin{equation}
F^{i a_{n-1}a_n} \varepsilon_{a_1\ldots a_n}\gamma^{a_1\ldots
a_{n-2}} (\mathbbm{1}-(-)^n w_\eta\gamma^*) \otimes e_i =0.
\label{nessuffD}
\end{equation}
Here the sum over the $a_*$ is over the unbarred indices with
respect to the complex structure given by the pure spinor $\eta$,
and the sum over $i$ is over the complete set of indices.
\end{proposition}

\begin{remark}\ \ {}
\begin{itemize}\itemsep=0pt
\item (\ref{nessuffD}) is solved by $F\in
(\mathbf{n}\oplus\bar{\mathbf{n}}) \otimes (
\mathfrak{su}(n)\oplus\Lambda^{0,2} )$. Of course, the strongly
torsion free Example~\ref{exsun} f\/its into this discussion.
\item Propositions \ref{proppurecurv} and \ref{proppuretor} give
the conditions on the connection such that the parallel pure
spinor yields a dif\/ferential.
\end{itemize}
\end{remark}

We will make a short comment on the twisted case. Consider  a
doubled spinor bundle. Suppose there are two pure spinors $\xi,
\hat\xi\in\Gamma S$, and let $\Xi=(\xi,\hat\xi)$ be one  parallel
spinor of the twisted bundle. Furthermore, suppose that the two
null-spaces def\/ined by $\xi$ and $\hat\xi$ intersect
transversally\footnote{This is true for the parallel pure spinor
and its charge conjugated counterpart  in the case of Levi-Civita
connection of holonomy $\mathfrak{su}(n)$. In this case
$\Xi=(\xi,\xi^C)$ is real and  $\mathfrak{B}(R;\Xi,\Xi)$ does not
vanish. This has been used in \cite{Klinker4} to show that the
real supersymmetric Killing structure  is not f\/inite in the case
of quaternionic spin representation where a twist of the spinor
bundle is necessary to yield a real structure. Nevertheless, it
has been shown that in this case there exist two isomorphic
f\/inite sub-structures.}. The necessary condition for
$\mathfrak{B}(R;\Xi,\Xi)$ to vanish is $R=0$. Now suppose that the
null spaces of the two spinors have non empty intersection $N$ and
the tangent space splits orthogonally into $T=N\oplus N^\perp$,
i.e.\ $\Lambda^2T=\Lambda^2N\oplus N\otimes N^\perp\oplus
\Lambda^2N^\perp $. Then the necessary condition reduces and only
the part of curvature which acts on $\Lambda^2N^\perp$  has to
vanish.


\subsection{Torsion freeness  from brane metrics}\label{examplebranes}

We consider a Lorentz\-ian manifold $M=(\r^D,g)$ such that the
coordinates are orthogonal with respect to the metric $g$.
Furthermore we consider a spinor connection $D^C$  which is
determined by a single $q$-form $F$. This $q$-form is  Hodge-dual
to a vector f\/ield $X$, where the Hodge-duality is with respect
to only one part  of the whole space. Furthermore the metric $g$
shall  depend on this vector f\/ield in such
 way that the Christof\/fel symbols obey
$\Gamma_{ABC}\propto X_Ag_{BC}$. We take $X$ to be the gradient of
a function $f$ and use the following ansatz for the metric on
$\r^D$:
\begin{equation}\label{pbranemetric}
g= f^2_\mu(x,y)\,(dx^\mu)^2+ f^2_i(x,y)\,(dy^i)^2,
\end{equation}
where $\big(x_\mu,y_m\big)_{0\leq \mu\leq p,1\leq m\leq d}$ is a
partition of coordinates into a $(p+1)$-dimensional space-time
determined by $\big(x_\mu\big)$ and a $d$-dimensional space
determined by $\big(y_m\big)_{1\leq m\leq d}$

We discuss two choices for the $q$-form $F$. Either $q=p+2$ with
\begin{gather}
F_{\mu_1\ldots\mu_{p+1} m}= \epsilon_{\mu_1\ldots\mu_{p+1}}\partial_m f(y)\nonumber\\ 
\intertext{or $q=d-1$ with} F_{m_1\ldots m_{d-1}}=
\epsilon_{m_1\ldots m_{d-1} m}\delta^{mn}\partial_n f(y),
\label{magneticF}
\end{gather}
where the function $f$  depends on $\{y_m\}$ only. We call $F$
electric or magnetic f\/ield strength in the f\/irst or second
case, respectively. This notation is due to the fact that the two
forms are connected via $F^{(p+2)}\propto *_{D}F^{(d-1)}$. Which
values for $p$ are possible to yield an admissible connection in
one of the two cases may be checked using
Theorem~\ref{admissibleForm} and its extension
Theorem~\ref{admissibleForm2}.

\begin{remark}
This metric together with the $q$-form for low dimensions is
considered in the discussion of $p$-brane solutions of
supergravity. E.g.\ in dimension $D=11$ we have a 5-brane with
magnetic four-form or a 2-brane with electric four-form. More
general p-branes may be obtained by using a non-f\/lat metric in
the space-time part (pp-waves or AdS) or in the space part (see
for example \cite{AliKumar} or \cite{Blau1} and references
therein).
\end{remark}

We specialize our discussion to the case where  the metric is
determined by two functions which depend on $\{y^i\}$ only:
\begin{equation} \label{pbranemetric1}
g= f^2_1(y) dx^2+ f^2_2(y)  dy^2\,.
\end{equation}
We refer to  the coordinate frame by unchecked indices and to the
orthonormal frame by checked indices. The two frames are connected
by  $e_{\check{\mu}}= f_1(y)^{-1}\partial_\mu$, $e_{\check{m}}=
f_2(y)^{-1}\partial_m$ and $e^{\check{\mu}}= f_1(y)dx^\mu$,
$e^{\check{m}}= f_2(y) dx^m$.

The Levi-Civita connection of (\ref{pbranemetric1}) is determined
by the Christof\/fel symbols
$\Gamma_{ABC}=\Gamma_{CBA}=\frac{1}{2}\big(\partial_A g_{BC}
+\partial_C g_{BA} -\partial_Bg_{AC}\big)$,
\begin{gather*}
\Gamma_{\mu \nu i }= -\Gamma_{\mu i\nu}\, =\, \partial_i(\ln f_1) g_{\mu\nu},\nonumber\\
\Gamma_{ijk}= \partial_i(\ln f_2) g_{jk}+\partial_k(\ln f_2) g_{ij} - \partial_j(\ln f_2) g_{ki},\\
\Gamma_{\mu \kappa \nu}= \Gamma_{\mu ij}\, =\, \Gamma_{i\mu j}
=0\nonumber ,
\end{gather*}
and given by
$\nabla_A=\partial_A+\tfrac{1}{4}\Gamma_{ABC}\gamma^{BC}$ with
\begin{gather*}
\nabla_\mu  = \partial_\mu+\tfrac{1}{2}\Gamma_{\mu \nu
i}\gamma^{\nu i}
         = \partial_\mu+\tfrac{1}{2} \partial_i(\ln f_1) f_1 f_2^{-1} \gamma_{\check \mu}{}^{\check \imath},\nonumber \\
\nabla_i    =
\partial_i+\tfrac{1}{4}\Gamma_{i\mu\nu}\gamma^{\mu\nu}+\tfrac{1}{4}\Gamma_{ijk}\gamma^{jk}
         = \partial_i+\tfrac{1}{2}\partial_j(\ln f_2) \gamma_{\check\imath}{}^{\check \jmath}.
\end{gather*}
The additional part $-\mathcal{A}^C$ of the spinor connection
$D^C=\nabla-\mathcal{A}^C$ is determined by the $q$-form~$F$ and
given by a linear combination of $F_{\mu A_1\ldots
A_{q-1}}\gamma^{A_1\ldots A_{q-1}}$ and $F^{A_1\ldots
A_q}\gamma_{\mu A_1\ldots A_q}$.

{\em The magnetic case. } We consider the $(d-1)$-form $F$ cf.\
(\ref{magneticF}) and calculate
\begin{gather*}
-\mathcal{A}^C_\mu
    = \alpha F_{\mu A_1\ldots A_{d-2}}\gamma^{A_1\ldots A_{d-2}}
                  +\beta F_{A_1\ldots A_{d-1}}\gamma_{\mu}{}^{A_1\ldots A_{d-1}}\\
\phantom{-\mathcal{A}^C_\mu}{}    = \beta \epsilon_{i_1\ldots i_{d-1}j}\delta^{jk}(\partial_k f)\gamma_{\mu}\gamma^{ i_1\ldots i_{d-1}} \\
\phantom{-\mathcal{A}^C_\mu}{}    = (-)^d\beta (d-1)! f_2^2 (\det g_d)^{-\frac{1}{2}} g^{jk}(\partial_j f)\gamma_{\mu k} \gamma^{[d]}\\
\phantom{-\mathcal{A}^C_\mu}{}    = (-)^d\beta (d-1)! (\partial_j f) f_1 f_2^2f_2^{-1}f_2^{-d} \gamma_{\check\mu}{}^{\check\jmath} \gamma^{[d]}\\
\phantom{-\mathcal{A}^C_\mu}{}    = (-)^d\beta (d-1)! (\partial_j
f) f_1 f_2^{1-d} \gamma_{\check\mu}{}^{\check\jmath} \gamma^{[d]}
\end{gather*}
as well as
\begin{gather*}
-\mathcal{A}^C_i
    =  \alpha F_{i A_1\ldots A_{d-2}}\gamma^{A_1\ldots A_{d-2}}
          +\beta F^{A_1\ldots A_{d-1}}\gamma_{i A_1\ldots A_{d-1}}\\
\phantom{-\mathcal{A}^C_i}{}    =  \alpha \epsilon_{i j_1\ldots
j_{d-2}k}\delta^{kj}(\partial_j f)\gamma^{j_1\ldots j_{d-2}}
          +\beta \epsilon_{j_1\ldots j_{d-1}k}\delta^{jk}(\partial_j f)\gamma_{i}{}^{j_1\ldots j_{d-1}}\\
\phantom{-\mathcal{A}^C_i}{}    =  (-)^{d-1}\big(\alpha (d-2)! g^{jk}f_2^2(\det g_d)^{-\frac{1}{2}}(\partial_j f)\gamma_{ik}\\
\phantom{-\mathcal{A}^C_i=}{}        +\beta
g^{jk}g_{ii'}f_2^2(\det g_d)^{\frac{1}{2}}(\partial_jf)
        \epsilon_{j_1\ldots j_{d-1}k}\epsilon^{i' j_1\ldots j_{d-1}}\big) \gamma^{[d]}\\
\phantom{-\mathcal{A}^C_i}{}     =  (-)^{d-1}\big(\alpha (d-2)!
(\partial_jf)f_2^2(\det
g_d)^{-\frac{1}{2}}\gamma_{\check\imath}{}^{\check\jmath}
          +\beta (d-1)!(\partial_i f)f_2^2(\det g_d)^{-\frac{1}{2}}\big) \gamma^{[d]}\\
\phantom{-\mathcal{A}^C_i}{}     =  (-)^{d-1}\big(\alpha (d-2)!
(\partial_jf)f_2^{2-d} \gamma_{\check\imath}{}^{\check\jmath}
          +\beta (d-1)!(\partial_i f)f_2^{2-d}  \gamma^{[d]}\big).
\end{gather*}

From now on we suppose that at least one of the two factors in the
brane ansatz is even dimensional. The matrix $\gamma^{[d]}
=\frac{1}{d!}\epsilon_{i_1\ldots i_d}\gamma^{i_1\ldots i_d}$ is
connected to the volume element of the space factor in $M$ and
obeys
\begin{equation*}
\gamma^{[d]}\gamma_j    =(-)^{d+1}\gamma_j\gamma^{[d]} \qquad
\text{and}\qquad \gamma^{[d]}\gamma_\mu
=(-)^{d}\gamma_\mu\gamma^{[d]}.
\end{equation*}
We choose $\varepsilon\in\{1,i\}$ such that
$(\varepsilon\gamma^{[d]})^2=\mathbbm{1}$. Then
$\Pi^\pm=\frac{1}{2}\big(\mathbbm{1}\pm
\varepsilon\gamma^{[d]}\big)$ is the projection on one half of the
spinor bundle. When we denote the spinor bundle of $M$, of its
$d$-dimensional factor, and of its ($p+1$)-dimensional factor by
$S_D$, $S_d$ and $S_{p+1}$, respectively, we have
\begin{equation*}
S_D^{\pm}=\Pi^{\pm}(S_D)=
\begin{cases}
S_d^{\pm}\otimes S_{p+1}&\text{if $d$ even},\\[1ex]
S_d\otimes S_{p+1}^\pm&\text{ if $d$ odd}.
\end{cases}
\end{equation*}

Furthermore we suppose  that $f_1$, $f_2$, $f$ and $\alpha$,
$\beta$ obey
\begin{gather}
X_i := \partial_i(\ln f_1)f_1 f_2^{-1}  =
 (-)^d\delta_1 \frac{2}{\varepsilon}\, \beta (d-1)!(\partial_i f)f_1f_2^{1-d},\nonumber \\
Y_i := \partial_i(\ln f_2)          =
-(-)^d\delta_2\frac{2}{\varepsilon}\,\alpha(d-2)!
(\partial_if)f_2^{2-d} \label{specialchoice}
\end{gather}
for some choice of signs $\delta_1,\delta_2\in\{\pm1\}$.

\begin{remark}
(\ref{specialchoice}) can be obtained by the ansatz
$f_\ell(y)=e^{\alpha_\ell u(y)}$ which yields the following system
for the constants $\alpha_\ell$:
\begin{equation*}
\alpha_1 = (-)^d\delta_1\frac{2}{\varepsilon}\,
\beta(d-1)!\alpha_3, \qquad \alpha_2=
-(-)^d\delta_2\frac{2}{\varepsilon}\,\alpha (d-2)!\alpha_3, \qquad
\alpha_3= (d-2)\alpha_2.
\end{equation*}

For $d=5$ we deal with a four-form $F$ which leads to an
admissible connection when we have $\Delta_0\Delta_1=-1$. Then  a
possible solution for $\delta_1=-\delta_2=-1$ and $\varepsilon=i$
is $\beta=-\frac{i}{288}$, $\alpha=\frac{8i}{288}$,
$\alpha_1=-\frac{1}{6}$, $\alpha_2=\frac{1}{3}$, and $\alpha_3=1$.
In dimension eleven this  is the supergravity M5-brane solution.
\end{remark}

With (\ref{specialchoice}) the connection $D^C$ is given by
\begin{equation}\label{specialchoiceconnection}
D^C_\mu=
\partial_\mu+X_i\gamma_{\check\mu}{}^{\check\imath}\Pi^\pm, \qquad
D^C_i =
\partial_i+Y_j\gamma_{\check\imath}{}^{\check\jmath}\Pi^\pm
        +\delta_2\varepsilon\frac{(d-1)\beta}{2 \alpha} Y_i \gamma^{[d]}  .
\end{equation}
The signs $\delta_*$  in (\ref{specialchoice}) determine which
projection is present. Nevertheless, the projections should be the
same in both terms.

\begin{proposition}\label{braneproposition}
$1.$ The holonomy of the connection
\eqref{specialchoiceconnection} is given by
\begin{equation*}
\mathfrak{hol}=     \mathfrak{so}(d)\ltimes
\begin{cases}
\displaystyle
     (p+1)\cdot 2^{\frac{d-1}{2}}\cdot S_d
        &\text{if $d$ odd, (i.e.\ $(p+1)$ even)},
\\[1.5ex]
\displaystyle
    (p+1)\cdot 2^{\frac{d}{2}-1} \cdot \tilde S_d
        &\text{if } d\equiv 0\,{\text{\rm mod}} \, 4,
\\[1ex]
\displaystyle
    (p+1)\cdot 2^{\frac{d}{2}}\cdot \tilde S_d
        &\text{if } d\equiv 2\,{\text{\rm mod}}\, 4,
\end{cases}
\end{equation*}
Here $\mathfrak{so}(d)\subset\mathfrak{sl}(S_D^\pm)$ and for $d$
even $\tilde S_d$ denotes the $2^{\frac{d}{2}-1}$-dimensional (not
specified) half spinor representation $S^\pm_d$.

$2.$ The torsion of the  connection $D$ -- the charge conjugated
of \eqref{specialchoiceconnection} -- is given by
\begin{gather}
\mathcal{T}_{\mu\nu}     =  \delta_1\varepsilon f_1 f_2 X\cdot \gamma_{\check\mu\check\nu} \gamma^{[d]},\nonumber \\
\mathcal{T}_{ \mu i}     =  -\delta_1\varepsilon f_2 X_i
\gamma_{\check\mu}\gamma^{[d]}
             =   \delta_2\varepsilon(d-1)\beta\alpha^{-1} f_1 Y_i \gamma_{\check\mu} \gamma^{[d]},  \label{branetor}\\
\mathcal{T}_{ij}         = \delta_2\varepsilon f_2 Y_k
\gamma^{\check k}{}_{\check\imath\check\jmath}
\gamma^{[d]}.\nonumber
\end{gather}
\end{proposition}

\begin{proof}
The bracket  $[D_\mu,D_\nu]$ vanishes due to
$\Pi^\pm\gamma_{\check\mu\check\imath}
=\gamma_{\check\mu\check\imath}\Pi^\mp$ and $\Pi^\mp\Pi^\pm=0$
whereas $[D_\mu,D_i]$  is given by
\begin{gather*}
 \big[\partial_\mu+X_j\gamma_{\check\mu}{}^{\check\jmath}\Pi^\pm,
\partial_i+Y_j\gamma_{\check\imath}{}^{\check\jmath}\Pi^\pm
        +\delta_2\varepsilon\frac{(d-1)\beta}{2\alpha} Y_i\gamma^{[d]} \big] \\
\qquad{}= \big[
X_j\gamma_{\check\mu}{}^{\check\jmath}\Pi^\pm,Y_k\gamma_{\check\imath}{}^{\check
k}\Pi^\pm    \big]
    + \Big[ X_j\gamma_{\check\mu}{}^{\check\jmath}\Pi^\pm, \delta_2\varepsilon
    \frac{(d-1)\beta}{2\alpha}(d-1) Y_i\gamma^{[d]}\Big]
    -\partial_iX_j\gamma_{\check\mu}{}^{\check\jmath}\Pi^\pm \\
\qquad {}=  X_jY_k \big[
\gamma_{\check\mu}{}^{\check\jmath}\Pi^\pm,\gamma_{\check\imath}{}^{\check
k}\Pi^\pm\big]
    +\delta_2\frac{(d-1)\beta}{2\alpha} Y_i X_j
       \big[\gamma_{\check\mu}{}^{\check\jmath}\Pi^\pm, \varepsilon\gamma^{[d]}\big]
    -\partial_iX_j\gamma_{\check\mu}{}^{\check\jmath}\Pi^\pm \\
\qquad =  X_jY_k \big( \gamma_{\check\mu}{}^{\check\jmath}\Pi^\pm
\gamma_{\check\imath}{}^{\check k}\Pi^\pm
   -    \gamma_{\check\imath}{}^{\check k}\Pi^\pm \gamma_{\check\mu}{}^{\check\jmath}\Pi^\pm\big) \\
\qquad\phantom{=}{} +\delta_2\frac{(d-1)\beta}{2\alpha} Y_i X_j
\big(\gamma_{\check\mu}{}^{\check\jmath}\Pi^\pm
\varepsilon\gamma^{[d]}
   - \varepsilon\gamma^{[d]} \gamma_{\check\mu}{}^{\check\jmath}\Pi^\pm \big)
      -\partial_iX_j\gamma_{\check\mu}{}^{\check\jmath}\Pi^\pm \\
\qquad {}= X_jY_k \gamma_{\check\mu}{}^{\check\jmath}
\gamma_{\check\imath}{}^{\check k}\Pi^\pm
    \pm \delta_2\frac{(d-1)\beta}{2\alpha} Y_i X_j  \gamma_{\check\mu}{}^{\check\jmath}\Pi^\pm
    -\partial_iX_j\gamma_{\check\mu}{}^{\check\jmath}\Pi^\pm \\
\qquad{}=  X_jY_k \gamma_{\check\mu}
             ( -\delta^{\check\jmath}_{\check\imath}\gamma^{\check k}+g^{\check j\check k}\gamma_{\check\imath} )\Pi^\pm
    \pm \delta_2\frac{(d-1)\beta}{2\alpha}Y_i X_j  \gamma_{\check\mu}{}^{\check\jmath}\Pi^\pm
    -\partial_iX_j\gamma_{\check\mu}{}^{\check\jmath}\Pi^\pm \\
\qquad{}=  -X_iY_j \gamma_{\check\mu}{}^{\check\jmath}\Pi^\pm+
X_kY^k f_2^2  \gamma_{\check\mu\check\imath}\Pi^\pm
     \pm \delta_2\frac{(d-1)\beta}{2\alpha} Y_i X_j  \gamma_{\check\mu}{}^{\check\jmath}\Pi^\pm
    -\partial_iX_j\gamma_{\check\mu}{}^{\check\jmath}\Pi^\pm \\
\qquad{}= \Big( \Big(
\frac{\pm\delta_2(d-1)\beta-2\alpha}{2\alpha}\, X_i Y_j  -
\partial_iX_j\Big)
     \gamma_{\check\mu}{}^{\check\jmath}
    + X_kY^k f_2^2  \gamma_{\check\mu\check\imath} \Big)\Pi^\pm.
\end{gather*}
Here we used  $\gamma_{\check\mu\check\imath}\gamma^{[d]}
=-\gamma^{[d]}\gamma_{\check\mu\check\imath}$ and $X_jY_k=X_kY_j$.
When we calculate $[D_\mu,D_i]$ we furthermore use
$\gamma_{\check\imath\check\jmath}\gamma^{[d]}=\gamma^{[d]}\gamma_{\check\imath\check\jmath}
$. This yields
\begin{gather*}
\Big[
\partial_i+Y_k\gamma_{\check\imath}{}^{\check k}\Pi^\pm
        +\delta_2\varepsilon\frac{(d-1)\beta}{2\alpha} Y_i\gamma^{[d]},
\partial_j+Y_\ell\gamma_{\check\jmath}{}^{\check\ell}\Pi^\pm
        +\delta_2\varepsilon\frac{(d-1)\beta}{2\alpha} Y_j\gamma^{[d]}
\Big] \\
\qquad{}=  (\partial_iY_\ell)
\gamma_{\check\jmath}{}^{\check\ell}\Pi^\pm
      -(\partial_jY_\ell) \gamma_{\check\imath}{}^{\check\ell}\Pi^\pm
     +\delta_2\varepsilon\frac{(d-1)\beta}{2\alpha}(\partial_iY_j-\partial_jY_i)\gamma^{[d]}
   +Y_kY_\ell \big[\gamma_{\check\imath}{}^{\check k}, \gamma_{\check\jmath}{}^{\check\ell} \big]\Pi^\pm \\
\qquad{}= 2 f_2^2 \partial_{[i}Y^\ell \gamma_{\check\jmath]
\check\ell} \Pi^\pm
     +2 Y_kY_\ell \big ( g_{\check\imath\check\jmath}\gamma^{\check k\check\ell}
        -\delta^{\check\ell}_{\check\imath}\gamma^{\check k}{}_{\check\jmath}
        +g^{\check k\check\ell}\gamma_{\check\imath\check\jmath}
        -\delta_{\check\jmath}^{\check k}\gamma_{\check\imath}{}^{\check\ell} \big) \Pi^\pm \\
\qquad{}= 2 \Big( f_2^2 \partial_{[i}Y^k \gamma_{\check\jmath]
\check k}
    + f_2^{-2}Y_kY^k \gamma_{\check\imath\check\jmath}
        +2Y_kY_{[i} \gamma_{\check\jmath]}{}^{\check k}
    \Big)\Pi^\pm .
\end{gather*}
We have two dif\/ferent families of generators for the holonomy
algebra: f\/irst
$\big\{\gamma_{\check\imath\check\jmath}\Pi^{\pm}\big\}$ and
second $\big\{\gamma_{\check\mu\check\jmath}\Pi^{\pm}\big\}$. The
f\/irst one generates a $\mathfrak{so}(d)$ sub algebra of
$\mathfrak{sl}(S^\pm_D)\subset \mathfrak{sl}(S_D)$.

Suppose $d$ is odd. The action of $\mathfrak{so}(d)$ on the second
family generates the commuting set
\begin{equation}\label{commuting}
{\rm
span}\big\{\gamma_{\check\mu\check\imath_1\ldots\check\imath_r}\Pi^{\pm}\,\big|\,
r \text{ odd}\big\} \simeq C\ell^{\text{odd}}_d.
\end{equation}
The action of $\mathfrak{so}(d)$ on this set is given by right
multiplication
\begin{equation*}
[\gamma_{\check\imath\check\jmath}\Pi^\pm,\gamma_{\check\mu\check
\imath_1\ldots\check\imath_r}\Pi^\pm] =
-\gamma_{\check\mu\check\imath_1
\ldots\check\imath_r}\gamma_{\check\imath\check\jmath}\Pi^\pm .
\end{equation*}
As a spin module via right (or left) multiplication the Clif\/ford
algebra is isomorphic to a direct sum of copies of the minimal
spinor representation $S_d$ and so is the $2^{d-1}$-dimensional
odd part due to Spin-invariance.

The minimal representation $S_d$ is of dimension
$2^{\frac{d-1}{2}}$. Therefore, the commuting set is isomorphic to
$(p+1)\cdot 2^{\frac{d-1}{2}} S_d$ as representation space.

Suppose $d$ is even. Consider once more the set generated by the
action of $\mathfrak{so}(d)$ on the second family. If
$(\gamma^{[d]})^2=\mathbbm{1}$  ($d\equiv0\,{\text{mod}}\,4$) we
have $\gamma^{\text{odd}}\Pi^{\pm}\propto
\gamma^{\text{odd}}\gamma^{[d]}\Pi^{\pm}=\gamma^{\text{odd}}\Pi^{\pm}$
in the other case ($d\equiv 2\,{\text{mod}}\,4$ )  there is an
extra $i$-factor in the proportionality. Therefore, the commuting
set is (\ref{commuting}) of dimension $2^{d-1}$ if
$d\equiv2\,{\text{mod}}\, 4$, and only one half of this if
$d\equiv 0\,{\text{mod}}\,4$ due to the duality above.

The minimal  representation $S_d^\pm$ is of  dimension
$2^{\frac{d}{2}-1}$. Therefore, as representation space the
commuting set is isomorphic to $(p+1)\cdot 2^{\frac{d}{2}-1}\tilde
S_d$ if $d\equiv 0\,{\text{mod}}\,4$  and to $(p+1)\cdot
2^{\frac{d}{2}}\tilde S_d$ if $d\equiv 2\,{\text{mod}}\,4$.

The torsion of the admissible connection $D=\nabla+\mathcal{A}$ is
given by
$\mathcal{T}_{AB}=ad^C_{\mathcal{A}_A}\gamma_B=\mathcal{A}_A\gamma_B+\gamma_B\mathcal{A}^C_A$.
We have
\begin{equation*}
-\mathcal{A}^C_\mu  =
\frac{\delta_1\varepsilon}{2}X_i\gamma_{\check\mu}{}^{\check\imath}\gamma^{[d]},\qquad
-\mathcal{A}^C_i     =
\frac{\delta_2\varepsilon}{2}Y_j\gamma_{\check\imath}{}^{\check\jmath}\gamma^{[d]}
                +\frac{\delta_2\varepsilon(d-1)\beta}{2\alpha}Y_i\gamma^{[d]}  .
\end{equation*}
Due to  Theorem \ref{admissibleForm} we have
$\Delta_{d-1}\Delta_1=-1$ or equivalently $\Delta_d\Delta_0=(-)^d$
which yields
 \begin{gather*}
-\mathcal{A}_\mu    =\
\frac{\delta_1\varepsilon}{2}X_i\big(\gamma_{\check\mu}{}^{\check\imath}\gamma^{[d]}\big)^C
            =   (-)^d\frac{\delta_1\varepsilon}{2}X_i \gamma_{\check\mu}{}^{\check\imath}\gamma^{[d]},\\
-\mathcal{A}_i      =  \frac{\delta_2\varepsilon}{2}Y_j
\big(\gamma_{\check\imath}{}^{\check\jmath}\gamma^{[d]}\big)^C
                +\frac{\delta_2\varepsilon(d-1)\beta}{2\alpha}Y_i (\gamma^{[d]})^C\\
  \phantom{-\mathcal{A}_i }{}          = -(-)^d \frac{\delta_2\varepsilon}{2}Y_j \gamma_{\check\imath}{}^{\check\jmath}\gamma^{[d]}
                +(-)^d\frac{\delta_2\varepsilon(d-1)\beta}{2\alpha}Y_i \gamma^{[d]}.
\end{gather*}
This is used to calculate the torsion of the brane connection:
\begin{gather*}
\mathcal{T}_{\mu\nu}
     = \mathcal{A}_\mu\gamma_\nu+\gamma_\nu\mathcal{A}^C_\mu
      =  -\frac{\delta_1\varepsilon}{2}X_i
        \big( (-)^d \gamma_{\check\mu}{}^{\check\imath}\gamma^{[d]}\gamma_\nu
            +\gamma_\nu \gamma_{\check\mu}{}^{\check\imath}\gamma^{[d]}\big) \\
\phantom{\mathcal{T}_{\mu\nu}}{}=
-\frac{\delta_1\varepsilon}{2}f_1 X_i
        \big( -\gamma_{\check\mu}\gamma_{\check\nu}\gamma^{\check\imath} \gamma^{[d]}
            +\gamma_{\check\nu} \gamma_{\check\mu}\gamma^{\check\imath}\gamma^{[d]}\big)
    =  \delta_1\varepsilon f_1f_2 X\cdot \gamma_{\check\mu\check\nu} \gamma^{[d]}
\end{gather*}
as well as
\begin{gather*}
\mathcal{T}_{\mu i}
    = \mathcal{A}_\mu\gamma_i+\gamma_i\mathcal{A}^C_\mu
 =  -\frac{\delta_1\varepsilon}{2}X_j
        \big( (-)^d \gamma_{\check\mu}{}^{\check\jmath}\gamma^{[d]}\gamma_i
            +\gamma_i \gamma_{\check\mu}{}^{\check\jmath}\gamma^{[d]}\big) \\
\phantom{\mathcal{T}_{\mu i}}{} = \frac{\delta_1\varepsilon}{2}f_2
X_j
        \big( \gamma^{\check\jmath}\gamma_{\check\imath} +\gamma_{\check\imath}\gamma^{\check\jmath}\big)
                         \gamma_{\check\mu}\gamma^{[d]}
 =  -\delta_1\varepsilon f_2 X_i \gamma_{\check\mu}\gamma^{[d]}.
\end{gather*}
Last but not least we have
\begin{gather*}
\mathcal{T}_{ij}
    = \mathcal{A}_i\gamma_j +\gamma_j \mathcal{A}^C_i \\
\phantom{\mathcal{T}_{ij} }{}    =
\frac{\delta_2\varepsilon}{2}Y_k
        \big( (-)^d\gamma_{\check\imath}{}^{\check k}\gamma^{[d]}\gamma_j
            -\gamma_j \gamma_{\check\imath}{}^{\check k}\gamma^{[d]}\big)
       - \frac{\delta_2\varepsilon(d-1)\beta}{2\alpha}Y_i
        \big( (-)^d \gamma^{[d]}\gamma_j +\gamma_j \gamma^{[d]}\big) \\
\phantom{\mathcal{T}_{ij} }{}= -\frac{\delta_2\varepsilon}{2}f_2
Y_k
        \big( \gamma_{\check\imath}{}^{\check k}\gamma_{\check\jmath}
            +\gamma_{\check\jmath} \gamma_{\check\imath}{}^{\check k}\big) \gamma^{[d]}
=  \delta_2\varepsilon f_2 Y_k  \gamma^{\check
k}{}_{\check\imath\check\jmath} \gamma^{[d]}.\tag*{\qed}
\end{gather*}\renewcommand{\qed}{}
\end{proof}

\begin{corollary}
The spinors which are parallel with respect to the connection
\eqref{specialchoiceconnection} form a subspace of the kernel of
$\Pi^\pm$. Explicitly we have $\eta(y)=f(y)\eta_0$ with constant
$\eta_0\in S_D^\mp $ and $f$ obeys $\partial_if =\pm
\delta_2\frac{(d-1)\beta}{2 \alpha} Y_i f $.
\end{corollary}

{\em The electric case.}   Due to the fact that the electric
$(p+2)$-form is dual to the magnetic $(d-1)$-form we will only
give a rough sketch of what is used to get a similar result. We
will assume that at least one of the factors is even dimensional.
Then we have the duality relation induced by
$\gamma^{[d]}\gamma^{[D]}\propto \gamma^{[p+1]}$.
 For a suitable choice of $X$ and $Y$ we get
\begin{equation*}
D^C_\mu = \partial_\mu
+X_i\gamma_{\check\mu}{}^{\check\imath}\hat\Pi^{\pm}, \qquad D^C_i
= \partial_\mu
+Y_j\gamma_{\check\mu}{}^{\check\jmath}\hat\Pi^{\pm}+\alpha
Y_i\gamma^{[p+1]}
\end{equation*}
which is of the same type as in the magnetic case. The projections
are given by
\[
\hat\Pi^\pm : S_D =S_d\otimes S_{p+1} \to S_D^\pm =
\begin{cases}
S_d\otimes S_{p+1}^\pm  &\text{ if ($p+1$) is even}, \\[1ex]
S_d^\pm \otimes S_{p+1} &\text{ if ($p+1$) is odd}.
\end{cases}
\]
The expressions for the holonomy and the torsion can be taken
directly from Proposition~\ref{braneproposition}.

In the remaining part of this section we analyze in what way we
have to restrict the set of parallel spinors to yield a torsion
free subset $\mathcal{K}$ in the sense of
Def\/inition~\ref{definitiontorsionfree}. I.e.\ we look for
solutions of $\mathfrak{D}(\mathcal{T},\eta,\xi)=0$ or
equivalently
\begin{gather}
 \mathcal{T}_{i\mu}\eta\wedge \gamma^\mu\xi +\mathcal{T}_{ij}\eta\wedge \gamma^j\xi +(\eta\leftrightarrow\xi) =0,
\nonumber\\
 \mathcal{T}_{\mu i}\eta\wedge \gamma^i\xi +\mathcal{T}_{\mu\nu}\eta\wedge \gamma^\nu \xi +(\eta\leftrightarrow\xi)=0
\label{torfreebrane1}
\end{gather}
with $\mathcal{T}_{AB}$ given by (\ref{branetor}). We discuss the
four summands separately and get:
\begin{gather*}
\mathcal{T}_{i\mu}\eta\wedge \gamma^\mu\xi
        = -\delta_1\varepsilon f_2X_i\gamma_{\check\mu}\gamma^{[d]}\eta\wedge \gamma^\mu\xi
        =   \delta_1\varepsilon f_2f_1^{-1}X_i \gamma_{\check\mu}\xi\wedge \gamma^{\check\mu}\gamma^{[d]}\eta \\
\phantom{\mathcal{T}_{i\mu}\eta\wedge \gamma^\mu\xi}{} =
\delta_1\varepsilon f_2 X_i
\gamma_{\check\mu}\gamma^{[d]}\xi\wedge \gamma^\mu\eta
        = - \mathcal{T}_{i\mu}\xi\wedge \gamma^\mu\eta,
\end{gather*}
where the second last equality holds because the spinors are of
the same chirality with respect to~$\gamma^{[d]}$.

From this result we see that the f\/irst summand in
(\ref{torfreebrane1}) vanishes after symmetrization over the
spinorial entries. With $Z_k=\delta_2\varepsilon f_2 Y_k $ we keep
on calculating
\begin{gather*}
\mathcal{T}_{ij}\eta\wedge \gamma^j\xi
  =   Z_k  \gamma^{\check k}{}_{\check\imath\check\jmath} \gamma^{[d]}\eta \wedge \gamma^j\xi
 =   Z_{\check k}\gamma_{\check\imath\check\jmath}\gamma^{\check k} \gamma^{[d]}\eta \wedge \gamma^{\check\jmath}\xi
        + Z_{\check k}\gamma_{[\check\imath} \delta^{\check k}_{\check\jmath]} \gamma^{[d]}\eta
        \wedge \gamma^{\check\jmath}\xi \\
\phantom{\mathcal{T}_{ij}\eta\wedge \gamma^j\xi }{} =
\gamma_{\check\imath\check\jmath}\, Z\cdot \gamma^{[d]}\eta \wedge
\gamma^{\check\jmath}\xi
    + \gamma_{\check\imath} \gamma^{[d]}\eta \wedge\, Z \xi
    - Z_{\check\imath}\, \gamma_{\check\jmath } \gamma^{[d]}\eta \wedge \gamma^{\check\jmath}\xi.
\end{gather*}
The last summand vanishes when we symmetrize with respect to
$\eta$ and $\xi$.  Furthermore we have
\begin{gather*}
T_{\mu i}\eta\wedge\gamma^i\xi
    = -\delta_1\varepsilon f_2 \gamma_{\check\mu}\gamma^{[d]}\eta\wedge (X\xi)
    = -(-)^d\delta_1\varepsilon f_2 \gamma^{[d]} \gamma_{\check\mu}\eta\wedge (X\xi)\\
\intertext{and} T_{\mu\nu}\eta\wedge\gamma^\nu\xi
    = \delta_1\varepsilon f_2 f_1 X\, \gamma_{\check\mu\check\nu} \gamma^{[d]} \eta  \wedge \gamma^{\nu}\xi
    = (-)^{d+1}\delta_1\varepsilon f_2  \gamma_{\check\mu\check\nu} \gamma^{[d]} (X\eta)  \wedge \gamma^{\check\nu}\xi .
\end{gather*}
If we put all this together and use $X\propto Z$ then  equations
(\ref{torfreebrane1}) reduce to
\begin{gather*}
  \gamma_{\check\imath\check\jmath}\gamma^{[d]} (X\eta) \wedge \gamma^{\check\jmath}\xi
    + \gamma^{[d]}\gamma_{\check\imath} \eta \wedge (X\xi) +(\eta\leftrightarrow\xi)  = 0,\nonumber\\
  \gamma_{\check\mu\check\nu} \gamma^{[d]} (X\eta) \wedge \gamma^{\check\nu}\xi
    + \gamma^{[d]} \gamma_{\check\mu}\eta\wedge (X\xi)
     +(\eta\leftrightarrow\xi) = 0.
\end{gather*}
We collect the brane-example in the following theorem.
\begin{theorem}
Consider the manifold $M$ which is $\r^{p+1}\times \r^d$  equipped
with the p-brane metric~\eqref{pbranemetric} and denote its spinor
bundle by $S$. Let $F$ be a magnetic $(d-1)$-form on $M$, i.e.\ it
is $*_d$-dual to a gradient field $X(y)$ on the transversal space
$\r^d$. The form $F$ and the metric are compatible such that they
define an admissible connection $D$ on the spinor bundle  cf.\
\eqref{specialchoiceconnection}. Then the space $\mathcal{K}$
given by
\[
\mathcal{K}=\big\{ \eta\in \Gamma S \,|\, D^C\eta=0, X\eta=0\big\}
\]
 is admissible and torsion free.
\end{theorem}

\section{Outlook}
As stated in the introduction admissible spinorial  connections,
i.e.\ connections with further symmetry condition on its torsion
c.f.\ Def\/inition~\ref{ctorsion}, are basic objects when we look
for inf\/initesimal automorphisms of the underlying manifold
constructed from parallel spinors, compare Theorem~\ref{KillingA}.
 This condition may be relaxed by considering admissible pairs cf.\ Def\/inition~\ref{defadmissible}.
 The notion of torsion enters naturally, when we look at commutators of vector
 f\/ields on supermanifolds constructed from the spinor bundle.
 This will be one tool in constructing a purely geometric
 representation of the supersymmetry algebra extending the work
 of \cite{AlCorDevSem} or \cite{Klinker4}. Work on this construction  is in progress.

\begin{appendix}

\section[Useful identities and symmetries for Clifford multiplication
 and charge conjugation]{Useful identities and symmetries for Clif\/ford multiplication\\
  and charge conjugation}\label{gammaappendix}

In this appendix we collect some identities concerning gamma
matrices as well as some properties of the symmetry of the
morphisms~(\ref{projectionC})\footnote{We note that most of the
formulas are valid without additional $(\det g)$-factors only if
the indices belong to an orthonormal frame (compare the
calculations in Section~\ref{examplebranes}).}.

For the Clif\/ford multiplication we use the convention
$\gamma_{\{\mu}\gamma_{\nu\}}=-g_{\mu\nu}$ which yields
\begin{equation}\label{iden}
\gamma_{\mu_1\ldots\mu_k}\gamma^{\nu_1\ldots\nu_\ell} =
\sum_{m=0}^{{\min}\{k,\ell\}}
\frac{(-)^{\frac{m(m-2k-1)}{2}}k!\ell!}{m!(k-m)!(\ell-m)!}
\delta^{[\nu_1\ldots \nu_m}_{[\mu_1\ldots\mu_m}
\gamma_{\mu_{m+1}\ldots\mu_k]}{}^{\nu_{m+1}\ldots\nu_\ell]}.
\end{equation}
We have
\begin{gather}
\gamma_{\mu_1\ldots \mu_k} =
\frac{1}{(D-k)!}(-)^{\frac{k(k+1)}{2}}(-)^{\frac{D(D+1)}{2}}
\epsilon_{\mu_1\ldots\mu_{D}}\gamma^{\mu_{k+1}\ldots\mu_{D}}\gamma^{[D]}
\label{duality}\\
\intertext{with} \gamma^{[D]}:= \gamma^1\cdots\gamma^{D}
    =\frac{1}{D!}\epsilon_{\mu_1\ldots\mu_{D}}\gamma^{\mu_1}\cdots\gamma^{\mu_{D}}.\nonumber
\end{gather}
This matrix obeys $(\gamma^{[D]})^2=(-)^{\frac{D(D+1)}{2}+t}$
 where $t$ denotes the amount of time-like directions in the metric.
For $D$ odd $\gamma^{[D]}$ is proportional to $\mathbbm{1}$. For
$D=2n$ even we def\/ine the modif\/ied matrix
\begin{equation*}
\gamma^*  =
\begin{cases}
\gamma^{[2n]} & \tilde\sigma\equiv0\,{\text{mod}}\,4, \\
i\, \gamma^{[2n]} & \tilde\sigma\equiv2\,{\text{mod}}\,4 ,
\end{cases}
\end{equation*}
where $\tilde\sigma$ denotes the signature of the metric $g$. It
obeys
\begin{equation*}
\gamma^*\gamma^{(k)}=(-)^k\gamma^{(k)}\gamma^*\qquad\text{and}\qquad
(\gamma^*)^2=\mathbbm{1}
\end{equation*}
and yields a splitting of the spinors in the two eigenspaces
$S=S^+\oplus S^-$.

The symmetry property (\ref{delta}) implies
\begin{gather*}
\Delta_k= -1          \   \Leftrightarrow\   k\in\{4m-\Delta_1, 4m+1+\Delta_0\}, \\
\Delta_0\Delta_k=  -1   \    \ \Leftrightarrow\  k\in\{4m+2, 4m-\Delta_0\Delta_1\},  \\
\Delta_1\Delta_k=-1 \ \Leftrightarrow\ \Delta_0\Delta_{k-1}= (-)^k
                \ \Leftrightarrow\   k\in\{4m+3, 4m+1+\Delta_0\Delta_1\}.
\end{gather*}
The symmetries $\Delta_k$ and $\Delta_{D-k}$ are connected via
\begin{equation*}
\Delta_{D-k}=(-)^{\frac{D(D-1)}{2}}(-)^{Dk}(-)^{k}(\Delta_0\Delta_1)^D\Delta_k.
\end{equation*}
This yields
\begin{alignat}{3}
& \Delta_{k} = (-)^{n+k}\Delta_{D-k} =
\Delta(\gamma^{(D-k)}\gamma^*)\quad  && \text{if }D=2n\text{
even},&
\label{starsymmetry}\\
& \Delta_{k} = \Delta_{D-k}             && \text{if }D=2n+1\text{
odd}.&\nonumber
\end{alignat}
Introducing the complex coordinates
$\gamma^a=\gamma^a+i\gamma^{a+n}$ and $\gamma^{\bar
a}=\gamma^a-i\gamma^{a+n}$,  for $a,\bar a=1,\ldots, n$, yields
\begin{gather*}
\gamma^{\{a}\gamma^{b\}}= \gamma^{\{\bar a}\gamma^{\bar b\}}=0,
\qquad
\gamma^{\{a}\gamma^{\bar b\}}=-2g^{a\bar b},  \\
(-)^{\frac{\sigma(\sigma-1)}{2}}\gamma^*=\gamma^{1\bar1}\cdots\gamma^{n\bar
n}
 =(\mathbbm{1}+\gamma^1\gamma^{\bar 1})\cdots (\mathbbm{1}+\gamma^n\gamma^{\bar n}),\\
\gamma^{1\ldots n}\gamma^*=\gamma^{1\ldots n},\qquad \gamma^{\bar
1\ldots\bar n}\gamma^*=(-)^n\gamma^{\bar 1\ldots\bar n}.
\end{gather*}

We use the following modif\/ied Pauli-matrices if we are forced to
modify the charge conjugation to change symmetries:
\begin{gather*}
\tau_0=\sigma_0=\begin{pmatrix}1&\\&1\end{pmatrix},\qquad
\tau_1=\sigma_1=\begin{pmatrix}&1\\1&\end{pmatrix},\\
\tau_2=i\sigma_2=\begin{pmatrix}&1\\-1&\end{pmatrix},\qquad
\tau_3=\sigma_3=\begin{pmatrix}1&\\&-1\end{pmatrix}.
\end{gather*}
To these matrices we associate two kind of  signs. The f\/irst
sign is $\varepsilon_{ik}$ which we get by permuting two of the
matrices, i.e.\
  $\tau_i\tau_k=\varepsilon_{ik}\tau_k\tau_i$,
  and the second is $\varepsilon_k$ which indicates the symmetry of $\tau_k$
\begin{equation*}
\varepsilon_{ik}=\begin{pmatrix}1&1&1&1\\ 1&1&-1&-1\\
1&-1&1&-1\\1&-1&-1&1\end{pmatrix},\qquad
\varepsilon_k=\begin{pmatrix}1\\1\\-1\\1 \end{pmatrix}.
\end{equation*}
\end{appendix}

\subsection*{Acknowledgements}
The author would like to thank Vicente Cort\'{e}s
 and Mario Listing for useful discussions on the topic of this text.

\LastPageEnding

\end{document}